\documentclass[a4paper,11pt]{article}
\usepackage{latexsym}
\addtocounter{section}{-1}

\input amssym.def
\input amssym

\usepackage{theorem}
\newtheorem{lemma}{Lemma}[section]
\newtheorem{proposition}{Proposition}[section]
\newtheorem{theorem}{Theorem}[section]
\newtheorem{corollary}{Corollary}[section]

\def\BC{\begin{center}}
\def\EC{\end{center}}
\def\BP{\nid {\it Proof.} }

\def\EP{\hspace*{\fill}$\Box$

\vx
}

\def\BT{\begin{theorem}}
\def\ET{\end{theorem}}
\def\BPR{\begin{proposition}}
\def\EPR{\end{proposition}}
\def\BL{\begin{lemma}}
\def\EL{\end{lemma}}
\def\BCO{\begin{corollary}}
\def\ECO{\end{corollary}}
\def\BPI{\BC\begin{picture}}
\def\EPI{\end{picture}\EC}
\def\BTA{\begin{tabular}}
\def\ETA{\end{tabular}}
\def\BIT{\begin{itemize}
\vx
}
\def\EIT{\end{itemize}}
\def\IT{
\mx
\item }
\def\BM{
\vx

\ \ \ $}

\def\EM{$

\vx
\nid }

\def\al{\alpha}
\def\be{\beta}

\def\ep{\varepsilon}
\def\ga{\gamma}

\def\KK{{\Bbb K}}
\def\LL{{\Bbb L}}

\def\vi{\varphi}

\def\1{\hspace{0.1ex}}
\def\2{\hspace{0.2ex}}
\def\3{\hspace{0.3ex}}
\def\4{\hspace{0.4ex}}
\def\5{\hspace{0.5ex}}
\def\6{\hspace{0.6ex}}
\def\7{\hspace{0.7ex}}
\def\8{\hspace{0.8ex}}
\def\9{\hspace{0.9ex}}

\def\l2{\hspace{-0.2ex}}
\def\mx{\vspace{-0.9ex}}
\def\vx{\vspace{1ex}}
\def\vh{\vspace{.5ex}}

\def\aa{\IT[{\rm (a) }]}
\def\bb{\IT[{\rm (b) }]}
\def\cc{\IT[{\rm (\hspace{0.1ex}c\hspace{0.1ex}) }]}
\def\dd{\IT[{\rm (d) }]}
\def\ee{\IT[{\rm (\hspace{0.1ex}e\hspace{0.1ex}) }]}

\def\A{{\cal A}}
\def\B{{\cal B}}
\def\C{\2{\cal C}\!\2}

\def\E{{\cal E}}

\def\G{{\cal G}}
\def\H{{\cal H}}
\def\I{{\cal I}}
\def\J{{\cal J}}

\def\L{{\cal L}}
\def\M{{\cal M}}

\def\S{{\cal S}}

\def\U{{\cal U}}

\def\alq{\forall\3 }
\def\all{\mbox{\hspace{0.5ex} for all \hspace{0.5ex}}}
\def\and{\mbox{\hspace{0.5ex} and \hspace{0.4ex}}}
\def\da{\downarrow\! }

\def\iff{\mathop {\hspace{0.4em}\Leftrightarrow\hspace{0.4em}}}
\def\imp{\mathop {\hspace{0.4em}\Rightarrow\hspace{0.4em}}}
\def\im{$\!\imp\!$}

\def\inc{\subseteq}
\def\Inf{\bigwedge}
\def\In{\mbox{$\mathop\bigcap$\1}}
\def\int{\2\cap\!\2}
\def\iso{\simeq}
\def\l{\hspace{-0.3ex}\mid\hspace{-0.3ex}}
\def\la{\leftarrow}

\def\map{\mathop{\1\rightarrow\1}}

\def\nid{\noindent}
\def\nI{\2{\!\s-\hspace{-1ex} I}\2}
\def\nJ{\2{\s-\hspace{-1.2ex} J}\3} 
\def\ol{\overline}

\def\ra{\rightarrow}

\def\s-{\2\setminus\!\2}

\def\Sup{\bigvee}
\def\ua{\uparrow\! }

\def\ve{\vee}
\def\we{\wedge}

\begin{document}
\title{Categories of Contexts}
\author{\sc Marcel Ern\'e}
\date{\small Department of Mathematics\\
University of Hannover\\
D-30167 Hannover, Germany}
\maketitle

\begin{quote}
\BC
{\bf Abstract}
\EC
\nid Morphisms between (formal) contexts are certain pairs of maps, one 
between objects and one between attributes of the contexts in question.
We study several classes of such morphisms and the connections between them.
Among other things, we show that the category ${\bf CLc}$ of complete 
lattices with complete homomorphisms is (up to a natural isomorphism) 
a full reflective subcategory of the category  of contexts 
with so-called conceptual morphisms; the reflector associates with each 
context its concept lattice. On the other hand, we obtain a dual adjunction
between ${\bf CLc}$ and the category of contexts with 
so-called concept continuous morphisms. Suitable restrictions of the 
adjoint functors yield a categorical equivalence and a duality between 
purified contexts and doubly based lattices, and in particular, between
reduced contexts and irreducibly bigenerated complete lattices. 
A central role is played by continuous maps between closure spaces and
by adjoint maps between complete lattices.

\vx


\end{quote}

\newpage

\section{Introduction}



Fundamental in Formal Concept Analysis is the interplay between so-called (formal) 
{\em contexts}, constituted by certain {\em incidence} relations, and the 
associated {\em concept lattices} introduced by Wille \cite{W1,W2}, 
hence lattice-theoretical tools in the spirit of Birkhoff \cite{BL}. It is 
therefore of primary interest to elucidate the passage between contexts and 
concept lattices -- and specifically, to investigate the relevant functors between the involved 
categories. 
Natural candidates on the lattice side are maps that preserve arbitrary 
joins, or meets, or both.
Often, morphisms between contexts will be {\em pairs} of maps, because
contexts always have two ground sets, one of ``objects'' and one of 
``attributes''. Since either of these sets carries a natural
closure system (that of ``extents'' and that of ``intents''), it is 
rather obvious that {\em continuity} will play a crucial role in that 
setting (see \cite{ELR} for a theory of lattice representations
for closure spaces). Continuity is also the defining condition for 
``scalings'' in measurement theory (see \cite{W1}).

\vx

In the present note, we are mainly interested in {\em complete 
homomorphisms} between concept lattices. Since they preserve both joins and
meets, we certainly have to take pairs of continuous maps 
between the underlying contexts -- but that is not enough, as observed in 
\cite{EDM}: one needs a certain link between the two involved mappings.
This leads us to two essentially different but equally important notions,
that of conceptual morphisms and that of concept continuous
morphisms: given contexts $(G,M,I)$ and $(H,N,J)$, a pair of mappings 
$\al : G\map H$ and $\be : M\map N$ is {\em conceptual} iff it preserves 
incidence (i.e. $g\2 I\2 m$ implies $\al (g)\2 J\2 \be (m)$) and an object 
$h$ has the attribute $n$ whenever each object $\be$-forced by $h$ has each
attribute $\al$-forced by $n$ (where $m$ is $\al${\em -forced} by $n$ if $m$ 
holds for each object whose $\al$-image has the attribute $n$, and 
$g$ is $\be${\em -forced} by $h$ if $g$ has each attribute whose $\be$-image
holds for $h$). On the other hand, $(\al ,\be )$ is {\em concept continuous}
iff it reflects incidence (i.e. $\al (j)\2 J\2 \be (m)$ implies $g\2 I\2 m$),
an attribute $n$ holds for $\al (g)$ whenever each $\be$-generalization
of $n$ holds for $g$, and $\be (m)$ holds for an object $h$ whenever $m$ 
holds for each  $\al$-specialization of $h$ (where $g$ is an 
$\al$-specialization of $h$ if $\al (g)$ shares all attributes of $h$, 
and $m$ is a $\be$-generalization of $n$ if $\be (m)$ holds for all objects
with attribute $n$). The category ${\bf CLc}$ of
complete lattices with complete homomorphisms turns out to be a full 
reflective subcategory of the category  of contexts with conceptual 
morphisms -- just by passing from contexts to their concept lattices. But 
there is also a {\em dual} adjunction between the category ${\bf CLc}$ and the 
category of contexts with concept continuous morphisms. Various results on
subcontexts and their concept lattices are immediate consequences.
Modifying the adjoint functors, we shall arrive at a categorical 
equivalence and a duality between purified contexts and so-called doubly 
based lattices -- in particular, between reduced contexts and 
irreducibly bigenerated complete lattices.

\newpage

\section{Categories of Complete Lattices}

In this preliminary section, we summarize some definitions and known facts about 
ordered sets, complete lattices and morphisms between them.

\vx

Given an arbitrary map $\vi : X\map Y$, we shall denote by $\vi [A]$ the 
image of $A\inc X$ under $\vi$, and by $\vi^-[B]$ the preimage of $B\inc Y$ under 
$\vi$. 

\vx

\nid Recall that a map $\vi : P\map Q$ between (partially) ordered sets is

\vh

-- {\em order preserving} or {\em isotone} if $\ x\leq y \imp \vi (x)\leq \vi (y)$,

\vh

-- {\em order reflecting} or {\em antitone} if $\ \vi (x)\leq \vi (y) \imp x\leq y$, 

\vh

-- an {\em order embedding} if $\ x\leq y \iff \vi (x)\leq \vi (y)$.

\vx

\nid Furthermore, $\vi$ is {\em join-dense} if each element in the codomain $Q$
is a join (supremum, least upper bound) of elements in the range of $\vi$,
or equivalently, if for $q\not\leq r$ in $Q$ there is a $p\in P$ with
$\vi (p)\leq q$ but $\vi (p)\not\leq r$. Caution: the composition of two
join-dense (isotone) maps need not be join-dense! \\
Meets and {\em meet-dense} maps are defined dually.

\vx

Of particular importance for our considerations are adjoint maps and 
functors (see \cite{DP}, \cite{EOV} or \cite{Com} for the order-theoretical and 
\cite{AHS} for the categorical part). Here we only recall the basic 
notions and facts.
A pair of maps $\vi : P\map Q$ and $\psi : Q\map P$ between ordered sets is 
{\em adjoint} if 
\BM
\vi (p)\leq q \iff p\leq \psi (q)
\EM
for all $p\in P$ and $q\in Q$. In that situation, $\vi$ is the {\em left}
or {\em lower adjoint} of $\psi$, which in turn is the {\em right} or 
{\em upper adjoint} of $\vi\3$. By antisymmetry of the order relations,
lower and upper adjoints determine each other uniquely; we write $\vi^*$ 
for the upper adjoint of $\vi$, and $\psi_*$ for the lower adjoint of 
$\psi$. 
It is helpful to know that a lower adjoint is injective iff its upper
adjoint is surjective, and {\em vice versa}. An injective lower
adjoint $\vi$ (upper adjoint $\psi$) is always an order embedding and 
satisfies $\vi^*\circ\vi = id$ ($\psi_*\circ\psi = id$). Note also that 
any join-dense join-preserving map and dually any meet-dense meet-preserving map 
is already surjective. 
Moreover, a map between posets is an isomorphism iff 
it has both an upper and a lower adjoint and these two adjoints coincide. 

\vx

\vx

Given subsets $A,B$ of a poset $P$, we denote by $A^{\ua}$ the collection
of all upper bounds of $A$ and by $B^{\da}$ that of all lower bounds of $B$.
In particular,
\BM
\da x = \{ x\}^{\da} = \{ p\in P : p\leq x\}\ $ and $\ 
\ua x = \{ x\}^{\ua} = \{ p\in P : p\geq x\}
\EM
are the {\em principal ideal} and the {\em principal filter} generated by $x\in P$,
respectively. More generally, for $A,B\inc P$,
\BM
A^{\ua\da} = \In \{ {\da x} : A\inc {\da x}\}
\EM
is the {\em lower cut} generated by $A$, and
\BM
B^{\da\ua} = \In \{ {\ua x} : B\inc {\ua x}\}
\EM
is the {\em upper cut} generated by $B$. The {\em cuts} in the sense of
MacNeille \cite{McN} (generalizing Dedekind's cuts of rational numbers) are
the pairs $(A,B)$ with $A= B^{\da}$ and $B= A^{\ua}$. Ordered by 
$(A,B)\leq (C,D) \!\iff\! A\inc B \!\iff\! D\inc C$, they form a complete
lattice, the {\em Dedekind-MacNeille completion}, which is
isomorphic to the closure system of lower cuts and dually isomorphic to the
closure system of upper cuts (cf. \cite{BL}, \cite{EOV}, \cite{EDM}).

\vx

A map between posets is a lower adjoint iff it is 
{\em residuated} (or {\em dually residual}), i.e. preimages of principal 
ideals are principal ideals, and it is an upper adjoint iff it is {\em residual}
(or {\em dually residuated}), i.e. preimages of principal filters are 
principal filters. Similarly, a map between posets is called {\em lower (upper) 
cut continuous} if preimages of lower (upper) cuts are again lower (upper)
cuts. From \cite{EDM}, we cite:

\BT
Generally, one has the following implications:

\vh

\BTA {lllll}
residuated &$\imp$& lower cut continuous &$\imp$& join preserving\\
residual &$\imp$& $\!$upper cut continuous &$\imp$& $\!$meet preserving
\ETA

\vh

\nid and for maps between complete lattices, the converse implications 
hold, too.
The completion by cuts yields a reflector from the category of posets with lower 
(upper) cut continuous maps to the full subcategory of complete lattices with join 
(meet) preserving maps.

\ET

We denote by ${\bf CLc}$ the category of complete lattices and {\em complete
homomorphisms}, i.e. maps preserving arbitrary joins and meets.
On the other hand, we have the category ${\bf CLc_*}$ of complete lattices
and {\em doubly residuated maps}, i.e. maps $\vi$ possessing an upper 
adjoint $\vi^*$ which again has an upper adoint $\vi^{**}$, 
and the category ${\bf CLc^*}$ of complete lattices and {\em doubly 
residual maps} $\psi$, having a lower adjoint $\psi_*$ which 
again has a lower adjoint $\psi_{**}$. Passing to upper adjoints, 
one obtains a dual isomorphism between the categories ${\bf CLc_*}$ 
and ${\bf CLc}$, but also one between ${\bf CLc}$ and ${\bf CLc^*}$. 
Composing both duality functors, one arrives at an isomorphism between the 
categories ${\bf CLc_*}$ and ${\bf CLc^*}$, sending any doubly residuated 
map $\vi$ to $\vi^{**}$, and in the opposite direction, any doubly residual
map $\psi$ to $\psi_{**}$. 
Summarizing the previous remarks, we note that under the above duality 
functors, the following pairs of categories of complete lattices are duals of each 
other:

\begin{picture}(0,0)
\put(-28,-65){

\BTA {|l|l||l|l|}
\hline
$\!$category$\!\!$& morphisms                  & dual           & morphisms\\
\hline
\hline
${\bf CLc}$ & complete homomorphisms & ${\bf CLc}_*$ & 
doubly residuated maps\\
             &                            & ${\bf CL^*c}$       & 
doubly residual maps\\
\hline
${\bf CLcd}$ & surjective (=$\5$dense) & ${\bf CLc_*e}$ & 
inj.$\5$doubly residuated maps\\
             & complete homomorphisms & ${\bf CLc^*e}$       & 
inj.$\5$doubly residual maps\\
\hline
${\bf CLce}$ & injective (=$\5$embedding) & ${\bf CLc_*d}$ & 
surj.$\5$doubly residuated maps$\!$\\
             & complete homomorphisms  & ${\bf CLcd}^*$       & 
surj.$\5$doubly residual maps\\
\hline
${\bf CLis}$ & isomorphisms & ${\bf CLis}$ & isomorphisms\\
\hline

\ETA 
}
\end{picture}

\newpage

\section{Closure Spaces and Continuous Maps}

Since contexts and their concept lattices are always intimately related with 
certain closure structures, a few preliminary remarks about closure spaces 
and their morphisms are in order before starting the morphism theory 
for contexts and concept lattices. For more background concerning the 
interaction between closure spaces and complete lattices, we refer to 
\cite{ELR}.

\vx

A {\em closure space} $X$ is a set together with a {\em closure system}
$\A (X)$, that is, a collection of subsets that is closed under arbitrary 
intersections. It is common use to denote the underlying set by the same 
letter as the space; thus, $X = \In \emptyset \in \A (X)$. For each subset
$A$ of $X$, there is a last member of $\A (X)$ containing $A$, denoted by 
$\ol{A}$ and called the {\em closure} of $A$. Clearly, $\A (X)$ is a complete 
lattice in which arbitrary meets coincide with intersections (but joins not
always with unions). There is a canonical map from $X$ to $\A (X)$,
\BM
\eta_X : X\map \A (X), \ x\mapsto \ol{\{ x\}} .
\EM
Let us recall several equivalent definitions of {\em continuity} for maps 
between closure spaces (see e.g. \cite{ELR}):

\BT
\label{cont}
For a map $\al$ between closure spaces $X$ and $Y$, the 
following conditions are equivalent:
\BIT
\aa Preimages of closed sets under $\al$ are closed.

\bb $\al [\ol{A}]\inc \ol{\al [A]}$ for all $A\inc X$.

\cc There are adjoint maps
$\al^{\ra}\! : \A (X) \map \A (Y)$ and $\al^{\la} : \A (Y) \map \A (X)$ 
with $\ \al^{\ra}\! \circ \eta_X = \eta_Y \circ \al.$

\dd There is a join-preserving $\al^{\ra}\!  : \A (X) \map \A (Y)$ 
with $\al^{\ra}(\ol{\{ x\}}) = \ol{\{\al (x)\} }$.
\EIT

\nid Moreover, these maps $\4\al ^{\ra}$ and $\al^{\la}$ are
uniquely determined:
\BM
\al^{\ra}(A)= \ol{\al [A]}, \ \al^{\la}(C)= \al^{-}[C].
\EM
\ET

In order to determine under what conditions the maps $\al^{\la}$
and $\al^{\ra}$ are injective or surjective, respectively, we say
$\al$ is

\vx

-- {\em (strictly) dense} if for each $B\inc Y$ there is some $A\inc X$ with
$\ol{B}=\ol{\al [A]}$,

\vh

-- {\em full} if $\al^-[\ol{\al [A]}]\inc \ol{A}$ for all $A\inc X$,
i.e. $\al (x)\in \ol{\al [A]}$ implies $x\in\ol{A}$.

\vx

\nid We shall omit the word ``strictly''. The following facts are easily 
checked: 

\BL
\label{insur} 
A continuous map $\al$ between closure spaces $X$ and $Y$ is

\vh

\BTA {lllll}
-- dense &iff&$\!\!\al^{\ra}$ is surjective&iff& 
$\hspace{-1ex}\al^{\la}$ is injective,\\
-- full &iff&$\!\!\al^{\ra}$ is injective&iff& 
$\hspace{-1ex}\al^{\la}$ is surjective.\\
\ETA

\vh

\nid Furthermore, $\al$ is full and continuous iff it is initial, i.e.
$\A (X) = \al^{\la}[\A (Y)]$.
\EL

\section{Morphisms Between Contexts}

A {\em (formal) context} is a triple $\KK = (G,M,I)$ where $I$ is some 
{\em `` incidence''} relation between elements of $G$ ({\em ``objects''}) 
and elements of $M$ ({\em ``marks''} or {\em ``attributes''}), i.e. 
$I\inc G\!\times\! M$. For $A\inc G$ and $B\inc M$, we put
\BM
A^{\ua} = A^{I} = \{ m\in\!\2 M : g\2 I\2 m \all g\in A\},$

\ \ \ $B^{\da} = B_{I} = \{ \6 g\in G\6 : g\2 I\2 m \all\! m\in B\} .
\EM
Instead of $\{ g\}^{\ua}$ and $\{ m\}^{\da}$, we shall write $g^{\ua}$ and
$m^{\da}\3$, respectively.\\
The {\em complementary} relation $G\!\times\! M\s- I$ will be denoted by 
$\nI$. A {\em (formal) concept} of the context $\KK$
is a pair $(A,B)$ with $A\inc G$, $B\inc M$, $A=B^{\da}$ (the 
{\em ``extent''}) and $A^{\ua} = B$ (the {\em ``intent''}). Ordered by
\BM
(A,B)\leq (C,D) \iff A\inc C \iff D\inc B\2 ,
\EM
the concepts form a complete lattice, the so-called {\em concept lattice}
$\B\2 \KK\2$. By passing to the first or second components, this lattice is 
isomorphic to the closure system $\E \KK$ of all extents and dually 
isomorphic to the closure system $\I\2 \KK$ of all intents. Thus, concept
lattices are the natural generalization of Dedekind-MacNeille completions, 
replacing order relations by arbitrary relations.
Notice that the corresponding closure operators are $^{\ua\da}$ and
$^{\da\ua}$, respectively: indeed, $A^{\ua\da}$ is the least extent 
containing $A\inc G$, and $B^{\da\ua}$ is the least intent containing 
$B\inc M$.

\vx

As mentioned in the introduction, we are mainly interested in categorical 
aspects of Formal Concept Analysis, that is, in the investigation of 
suitable morphisms between contexts and the associated concept lattices. 
Naturally, context morphisms have to be 
certain pairs of maps, one between the objects and the other between the 
attributes. The choice of morphisms is not evident and may 
depend heavily on the intended investigations and results. From the 
Galois-theoretical point of view, it would be natural to consider pairs of 
maps with opposite directions. However, we shall not pursue that trace in 
the present note but focus on situations where both maps run into the same 
direction - an approach that leads to quite satisfactory results as well. 

\vx

Given two contexts $\KK = (G,M,I)$ and $\LL = (H,N,J)$, a pair 
$(\al, \beta )$ of maps $\al : G\map H$ and $\beta : M\map N$ will be 
referred to as a {\em mapping pair} or {\em (weak) concept morphism}. 
Let us list a few natural conditions on such mappings. In accordance with 
the corresponding general closure-theoretical definitions, we say $\al$ is 

\vh

-- {\em (extent) continuous} if preimages
of extents under $\al$ are extents,

\vh 

-- {\em (extent) dense} if for all $C\inc H$, there is an $A\inc G$ with
$C^{\ua\da} = \al [A]^{\ua\da}$,

\vh

-- {\em (extent) full} if for all $A\inc G$, 
$\al (g)\in \al [A]^{\ua\da}$ implies $g\in A^{\ua\da}$.

\vx

\nid Dually, $\beta$ is said to be  

\vh

-- {\em (intent) continuous} if preimages of intents under $\beta$ are 
intents,

\vh

-- {\em (intent) dense} if for all $D\inc N$, there is a $B\inc M$ with
$D^{\da\ua} = \be [B]^{\da\ua}$,

\vh

-- {\em (intent) full} if for all $B\inc M$, 
$\be (m)\in \be [B]^{\da\ua}$ implies $m\in B^{\da\ua}$.

\vx

\nid Extent continuous maps are often interpreted as {\em scalings} in 
the theory of measurement, in particular if the objects
of the codomain are numbers or numerical functions (see e.g. \cite{W1}).

\vx

Although every closure space $X$ may be regarded as an extent space, 
namely of the context $(X,\A(X),\in )$, there is a crucial difference between 
arbitrary closure spaces and extent or intent spaces: in the latter situation 
the various types of morphisms admit descriptions in first order terms, 
involving quantification over objects and attributes only, but not over subsets 
(like extents or intents). This reduction of complexity is one of the prominent 
advantages of Formal Concept Analysis (where contexts are regarded as 
``logarithms'' of their concept lattices). Note that statements like
\BM
h^{\ua}\inc\al (g)^{\ua}$ \ or \ $\al^{-}[n^{\da}\2]\inc m^{\da}
\EM
are expressible in first order terms (the former meaning that $h\2 J \2 n$
implies $\al (g) \2 J \2 n$, and the latter that $\al (g)\2 J\2 n$ implies
$g \2 I\2 m$).

\vx

In the subsequent lemmas, $(\al ,\be )$ always denotes a mapping pair 
between contexts $\KK = (G,M,I)$ and $\LL = (H,N,J)$. 

\BL
The following are equivalent:
\BIT
\mx
\aa $\al$ is extent continuous.

\bb $\al^- [n^{\da}\2 ]$ is an extent for each $n\in N$.

\cc $\al (g)\nJ n$ implies $g\nI m$ 
for some $m\in M$ with $\al^{-}[n^{\da}\2]\inc m^{\da}$.

\dd $\al [A^{\ua\da}]^{\ua} = \al [A]^{\ua}$ for all $A\inc G$.

\ee $\al [A^{\ua\da}]\inc \al [A]^{\ua\da}$ for all $A\inc G$.
\EIT
Dual characterizations hold for intent continuous maps.
\EL

\BP
(a)\im(b)\im(c). $\al (g)\nJ n$ means 
$g\not\in \al^- [n^{\da}] = \al^- [n^{\da}]^{\ua\da}$
(as $\al^- [n^{\da}]$ is an extent). Thus, there is an $m\in M$ with
$g\nI m$ but $\al^{-}[n^{\da}\2]\inc m^{\da}$.

\vh

\nid (c)\im(d). If $n\not\in \al [A^{\ua\da}]^{\ua}$ then 
$\al (g)\nJ n$ for some $g\in A^{\ua\da}$, hence
$g\nI m$ but $\al^{-}[n^{\da}\2]\inc m^{\da}$ for some
$m\in M$. It follows that $A\not\inc m^{\da}$, {\em a fortiori} 
$A\not\inc \al^{-}[n^{\da}\2]$, which means $n\not\in \al [A]^{\ua}$. By
contraposition, we obtain $\al [A]^{\ua} \inc \al [A^{\ua\da}]^{\ua}$, and 
the reverse inclusion is a consequence of $A\inc A^{\ua\da}$.

\vh

\nid The implication (d)\im(e) is clear. For (e)\im(a), see Theorem \ref{cont}.
\EP

Note that condition (e) may be reformulated as an ``implication between
implications'' (cf. \cite{GWW}); indeed, writing 
$A\ra B$ for $A^{\ua}\inc B^{\ua}$, i.e. $B\inc A^{\ua\da}$
(meaning that the objects of $B$ share all common properties of objects in $A$), 
we see that (e) is equi\-valent to
\BM
{\rm (d')} \ A\ra B$ implies $\al [A] \ra \al [B]\2 .
\EM 

\vx


\BL
\label{dense}
The map $\al$ is extent dense 

\vh

iff for each $h\in H$ there is a set $A\inc G$ with $h^{\ua}=\al [A]^{\ua}$

\vh

iff $\ h\nJ n$ implies 
$\al (g)\nJ n$ 
for some $g\in G$ with $ h^{\ua} \inc \al (g)^{\ua}  $.

\vx

\nid Dually, the map $\be$ is intent dense 

\vh

iff for each $n\in N$ there is a set $B\inc M$ with $n^{\da}=\be [B]^{\da}$

\vh

iff $\ h\nJ n$ implies 
$h\nJ \be (m)$
for some $m\in M$ with $ n^{\da} \inc \be (m)^{\da}  $.

\vx

\nid Consequently, both $\al$ and $\be$ are dense

iff for $h\nJ n$ there exist $g,m$ with
$\al (g)\nJ \be (m)$, $h^{\ua}\inc \al (g)^{\ua}$,
$n^{\da}\inc \be (m)^{\da}$.

\EL

\BP
If $\al$ is dense then for each $h\in H$ there is some $A\inc G$ such that
$h^{\ua\da} = \al [A]^{\ua\da}$, hence $h^{\ua} = \al [A]^{\ua}$. Assuming 
the latter equation and $h\nJ n$, 
we find a $g\in A$ with $\al (g)\nJ n$, whereas
$ h^{\ua} \inc \al (g)^{\ua} $.

\vh
 
On the other hand, assume that $h\nJ n$ implies 
$\al (g)\nJ n$ 
for some $g$ with $ h^{\ua} \inc \al (g)^{\ua}$.
In order to prove density of $\al$, consider an extent $C=C^{\ua\da}$ and
put $A=\al^- [C]$. We claim that $\al [A]^{\ua\da}=C$. If $h\not\in \al [A]^{\ua\da}$, choose $n\in N$ with 
$h\nJ n$ but $\al [A]\inc n^{\da}$, and then 
$g\in G$ with $\al (g)\nJ n$ but 
$h^{\ua}\inc\al (g)^{\ua}$.
Then we have $A\inc \al^- [n^{\da}]$ but $g\not\in A$ (otherwise 
$\al (g)\in \al [A]\inc n^{\da}$, i.e. $\al (g)\2 J\3 n$), hence 
$\al (g)\not\in C$. The assumption $h\in C$ leads to 
$C^{\ua}\inc h^{\ua}\inc\al (g)^{\ua}$ and $\al (g)\in C^{\ua\da}=C$, a
contradiction. By contraposition, $C\inc \al [A]^{\ua\da}\inc C$.
\EP

Note also that the map $\al$ is extent dense iff the set
$\{ \al (g)^{\ua\da} : g\in G\}$ is join-dense in the extent lattice
$\E \LL$, and dually, $\be$ is intent dense iff the set
$\{ \be (m)^{\da\ua} : m\in M\}$ is meet-dense in the intent lattice
$\I\2 \LL$.

\BL
\label{full}
The map $\al$ is extent full 

\vh

iff $\ g\nI m$ implies $\al (g)\nJ n$ 
for some $n\in N$ with $\al [m^{\da}\2 ]\inc n^{\da}$

\vh

iff $\al [A]\!\imp\! \al [B]$ entails $A\!\imp\! B$

\vh

iff each extent of $\KK$ is the preimage of an extent of $\LL$ under $\al$.

\vx

\nid Dual characterizations hold for intent fullness.




\EL

\BP
If $\al$ is full then $g\nI m$ implies $g\not\in
m^{\da}=m^{\da\ua\da}$ and so $\al (g)\not\in \al [m^{\da}]^{\ua\da}$, i.e.
$\al (g)\nJ n$ 
for some $n\in N$ with $\al [m^{\da}]\inc n^{\da}$.
Conversely, if the latter holds then for $g\in G$ and $A\inc G$ with 
$g\not\in A^{\ua\da}$, we find an $m\in A^{\ua}$ with 
$g\nI m$ and then an $n\in \al [m^{\da}]^{\ua}$ such
that $\al (g)\nJ n$, whence $\al (g)\not\in n^{\da}$ 
and, {\em a fortiori}, $\al (g)\not\in \al [A]^{\ua\da}$ (because $A\inc
m^{\da}$ and therefore $\al [A]^{\ua\da}\inc \al [m^{\da}]^{\ua\da}\inc 
n^{\da}$).


\nid The last two characterizations of fullness are straightforward.
\EP

\BCO
$\al$ is initial, i.e. extent continuous and full 

\vh

iff $A\ra B$ is equivalent to $\al [A]\map \al [B]$

\vh

iff the extents of $\KK$ are precisely the preimages of extents of $\LL\2$.

\ECO


\vx

We come now to the crucial definitions, relating both partners of a mapping pair
to each other. The mapping pair $(\al , \be )$ is called

\vh

-- {\em incidence preserving} if $g \2 I \2 m$ implies $\al (g) \2 J \3 \be (m)$,

\vh

-- {\em incidence reflecting} \2 if $\2 \al (g) \2 J \3 \be (m)$ implies 
$g \2 I \2 m$,

\vh

-- a {\em context embedding}\2 if it preserves and reflects incidence.

\vx

\nid It would be a bit more precise to speak of 
{\em quasi-embeddings} and to reserve the name {\em embeddings} to the case
where both mappings are one-to-one.
However, injectivity plays a minor role in the present study. 

\BL
\label{princid}
The mapping pair $(\al ,\be )$ preserves incidence 

\vh

iff $\ \al [A\1 ]^{\ua\da}\inc \be [A^{\ua}\3 ]^{\da}$ for all 
$A\inc G$

\vh

iff $\ \be [B]^{\da\ua}\inc \al [B^{\da}\2 ]^{\ua}$ for all $B\inc M$.
\EL

\BP
If $(\al ,\be )$ preserves incidence then 
$\be [A^{\ua}\2 ]\!\inc\! \al [A]^{\ua}$, 
hence $\al [A]^{\ua\da}\!\inc\! \be [A^{\ua}\2 ]^{\da}$.\\
Conversely, that inclusion implies $\be [g^{\ua}]\inc \al (g)^{\ua}$ 
for all $g\in G$,
which means that $(\al ,\be )$ preserves incidence. 
The other equivalence is shown dually.
\EP

\BL
\label{refinc}
The mapping pair $(\al ,\be )$ reflects incidence 

\vh

iff $\ \al^- [C^{\ua\da}\2 ]\inc \be^- [C^{\1\ua\1}\3 ]^{\da}$ for all 
$C\inc H$

\vh

iff $\ \be^- [D^{\da\ua}\2 ]\inc \al^- [D^{\da}\2 ]^{\ua}$ for all $D\inc N$.
\EL

\BP
For $g\in G\s- \be^- [C^{\1\ua}\2 ]^{\da}$ we find an $m\in M$ with 
$\be (m)\in\! C^{\ua}$ and $\ g\!\5\nI m$; if 
$(\al ,\be )$ reflects incidence then $\ g\nI m$ 
implies $\al (g)\nJ \be (m)$, whence 
$\al (g)\not\in C^{\ua\da}$. By contraposition, we get
$\al^- [C^{\ua\da}\2 ]\inc \be^- [C^{\ua}\2 ]^{\da}$. Conversely, if that
inclusion is fulfilled then $\al (g) \2 J \2 \be (m)$ implies 
$\al^- [\al (g)^{\ua\da}\2 ]\inc \be^- [\al (g)^{\ua}\2 ]^{\da}\inc m^{\da}$
(since $m\in \be^- [\al (g)^{\ua}]$),
and it follows that $g\in \al^- [\al (g)^{\ua\da}\2 ] \inc m^{\da}$, i.e.
$g\2 I\2 m$.
\EP

\vx

Even more important than the above properties of mapping pairs are certain
strong kinds of continuity. We say a mapping pair $(\al ,\be )$ is

\vh

-- {\em separately continuous} if both $\al$ and $\be$ are continuous,

\vh

-- {\em concept preserving} if $(A,B)\in \B\2 \KK$ implies 
$(\be [B]^{\da},\al [A]^{\ua} )\in \B\2 \LL$,

\vh

-- {\em conceptual} if it is separately continuous and concept preserving,

\vh

-- {\em concept continuous} if $\ (C,D)\in\B\2 \LL$ implies
$(\al^- [C],\be^- [D])\in\B\2 \KK$,

\vh

-- a {\em dense context embedding} if it is an embedding and $\al ,\be $ are
dense,

\vh

-- a {\em context isomorphism} if it is an embedding and
$\al ,\be$ are bijective.



The next result has been shown in \cite{EDM} for the case of order relations:

\BL
\label{conceptual}
The mapping pair $(\al ,\be)$ is conceptual 

\vh

iff $\ \al [A]^{\ua\da} = \be [A^{\ua}\2 ]^{\da}$ for all $A\inc G$ and
$\be [B]^{\da\ua} = \al [B^{\da}\2 ]^{\ua}$ for all $B\inc M$

\vh

iff $(\al ,\be )$ preserves incidence, 
$\be [\al^-[n^{\da}\2 ]^{\ua}\2 ]^{\da}\inc n^{\da}$ and
$\al [\be^-[h^{\ua}\2 ]^{\da}\2 ]^{\ua}\inc h^{\ua}$ 
\vh 

iff $(\al ,\be )$ preserves incidence and for $h\nJ n$, there are 
$\ g\nI m$ with 

\ \ \ \ $\al^- [n^{\da}\2 ]\inc m^{\da}$ and $\be^- [h^{\ua}\2 ]\inc g^{\ua}$.
\EL

\BP
Suppose $(\al ,\be )$ is conceptual. For $A\inc G$, the pair $(A^{\ua\da},A^{\ua})$
is a concept; hence $(\be [A^{\ua}\2 ]^{\da}, \al [A^{\ua\da}\2 ]^{\ua})$
is a concept, too. Thus $\be [A^{\ua}\2 ]^{\da} = 
\al [A^{\ua\da}\2 ]^{\ua\da} =\al [A]^{\ua\da}$ (by continuity of $\al$). 
The second equation is obtained analogously.

\vh

Now, if we assume the equations $\ \al [A]^{\ua\da} = \be [A^{\ua}\2 ]^{\da}$ and
$\be [B]^{\da\ua} = \al [B^{\da}\2 ]^{\ua}$ then by Lemma \ref{princid}, 
$(\al ,\be )$ preserves incidence; furthermore, we have\\
$\be [\al^-[n^{\da}\2 ]^{\ua}\2 ]^{\da} = 
\al [\al^- [n^{\da}\2 ]]^{\ua\da} \inc n^{\da\ua\da} = n^{\da}\ $ and
$\ \al [\be^-[h^{\ua}\2 ]^{\da}\2 ]^{\ua} =
\be [\be^- [h^{\ua}\2 ]]^{\da\ua} \inc h^{\ua}$.

\vh
 
On the other hand, if $\be [\al^-[n^{\da}\2 ]^{\ua}\2 ]^{\da}$ is contained in
$n^{\da}$ then $h\nJ n$ implies $h\not\in \be [\al^-[n^{\da}\2 ]^{\ua}\2 ]^{\da}$,
so we find an $m\in \al^- [n^{\da}\2 ]^{\ua}$ with $\be (m)\not\in h^{\ua}$;
it follows that $\al^- [n^{\da}\2 ]\inc m^{\da}$ and 
$\be (m)\not\in \al [\be^-[h^{\ua}\2 ]^{\da}\2 ]^{\ua}$; therefore, we find
a $g\in \be^- [h^{\ua}\2 ]^{\da}$ with $\al (g) \nJ \be (m)$, hence $g\nI m$
(by incidence preservation). 


\vh

Finally, let us suppose that $(\al ,\be )$ preserves incidence and for 
$h\nJ n$, there are 
$\ g\nI m$ with $\al^- [n^{\da}\2 ]\inc m^{\da}$ and 
$\be^- [h^{\ua}\2 ]\inc g^{\ua}$. Again by Lemma \ref{princid}, we have
$\al [A]^{\ua\da}\inc \be [A^{\ua}\2 ]^{\da}$ and 
$\be [B]^{\da\ua}\inc \al [B^{\da}\2 ]^{\ua}$.
Assume $h\not\in \al [A]^{\ua\da}\2$; then there is some $n\in \al [A]^{\ua}$
with $h\nJ n$. Choose 
$\ g\nI m$ as above. Then $m\not\in g^{\ua}$, and the
inclusion $\be^- [h^{\ua}\2 ]\inc g^{\ua}$ yields  
$h\nJ \be (m)$. But $n\in \al [A]^{\ua}$ means
$\al [A]\inc n^{\da}$, i.e. $A\inc \al^- [n^{\da}\2 ]\inc m^{\da}$,
whence $m\in A^{\ua}$ and so $\be (m)\in \be [A^{\ua}\2 ]$.
This together with $h\nJ \be (m)$ gives
$h\!\not\in\! \be [A^{\ua}\2 ]^{\da}$, proving the equality  
$\ \al [A]^{\ua\da} = \be [A^{\ua}\2 ]^{\da}$.
\EP 

\BL
\label{concont}
The mapping pair $(\al ,\be)$ is concept continuous 

\vh

iff $\al^- [C^{\ua\da}\2 ] = \be^- [C^{\ua}\2 ]^{\da}$ for all $C\!\inc\! H$ 
and $\be^- [D^{\da\ua}\2 ] = \al^- [D^{\da}\2 ]^{\ua}$ for all $D\!\inc\! N$

\vh

iff $\al^- [n^{\da}\2 ] = \be^- [n^{\da\ua}\2]^{\da}$ for all $n\in N$ and
$\be^- [h^{\ua}\2 ] = \al^- [h^{\ua\da}\2 ]^{\ua}$ for all $h\in H$

\vh

iff $(\al ,\be )$ reflects incidence, 
$\al [\be^- [n^{\da\ua}\2 ]^{\da}\2 ]\inc n^{\da }$ and
$\be [\al^- [h^{\ua\da}\2 ]^{\ua}\2 ]\inc h^{\ua }$ , i.e.  

\ \ \ \ for $\al (g) \nJ n$, there is an $m\in M$ with 
$g\nI m$ and $n^{\da}\inc \be (m)^{\da}\2$, and 

\ \ \ \ for $h\nJ \be (m)$, there is a $g\in G$ with 
$g\nI m$ and $h^{\ua}\inc \al (g)^{\ua}$.
\EL

\BP
If $(\al ,\be )$ is concept continuous then, observing that for arbitrary
$C\inc H$ the pair $(C^{\ua\da},C^{\ua})$ is a concept, we infer that
$(\al^- [C^{\ua\da}\2 ],\be^- [C^{\ua}\2 ])$ is a concept, too.
Hence $\al^- [C^{\ua\da}\2 ] = \be^- [C^{\ua}\2 ]^{\da}$, and dually for
$D\inc N$, $\be^- [D^{\da\ua}\2 ] =\al^- [D^{\da}\2 ]^{\ua}$.
Of course, the latter two equations entail 
$\al^- [n^{\da}\2 ] = \be^- [n^{\da\ua}\2]^{\da}$ for all $n\in N$ and
$\be^- [h^{\ua}\2 ] = \al^- [h^{\ua\da}\2 ]^{\ua}$ for all $h\in H$
(take $C=n^{\da}$ and $D=h^{\ua}$). 

\vh

Assume in turn the validity of these equations. 
If $h\nJ \be (m)$ then $m$ is not a member of $\be^- [h^{\ua}\2 ]=
\al^- [h^{\ua\da}\2 ]^{\ua}$. Thus, there exists a $g\!\in\! G$ with
$g\nI m$ and $\al (g)\in h^{\ua\da}$, i.e.
$h^{\ua}\inc \al (g)^{\ua}$. Dually, we find for 
$\al (g)\nJ n$ an $m\!\in\! M$ with 
$g\nI m$ and $n^{\da}\inc \be (m)^{\da}$.
As in Lemma \ref{refinc} we see that $(\al ,\be )$ reflects incidence.

\vh

Finally, suppose that the latter three conditions are fulfilled. 
In order to show that for any concept $(C,D)$ of $(H,N,J)$, the ``inverse
image'' $(\al^- [C^{\ua\da}\2 ],\be^- [D^{\da\ua}\2 ])$ is a concept of
$(G,M,I)$, we have to verify the equations 
$\al^- [D^{\da}\2 ]=\be^- [D]^{\da}$ and $\be^-[C^{\ua}\2 ]=\al^- [C]^{\ua}$.
If $\al (g)\not\in D^{\da}$ then $\al (g)\nJ n$ 
for some $n\in D$. By hypothesis, there is an $m\in M$ with  
$g\nI m$ and $n^{\da}\inc \be (m)^{\da}$, whence
$m\in \be^- [n^{\ua\da}\2 ]\inc \be^- [D^{\ua\da}\2 ] = \be^- [D]$
and therefore $g\not\in \be^- [D]^{\da}$. This proves the inclusion
$\be^- [D]^{\da}\inc \al^- [D^{\da}\2 ]$, and the other inclusion
$\al^- [D^{\da}\2 ] = \al^- [C^{\ua\da}\2 ] \inc \be^- [C^{\ua}\2 ]^{\da}=
\be^- [D]^{\da}$ follows from Lemma \ref{refinc}.
\EP

The last two lemmas confirm our characterizations of conceptual and 
concept continuous pairs from the introduction: the pair $(\al ,\be )$ is 
conceptual iff it preserves incidence and an object $h$ has the attribute 
$n$ whenever each object $g$ that is $\be$-forced by $h$ (i.e. 
$\be^- [h^{\ua}\2 ]\inc g^{\ua}$) has each attribute $m$ that is 
$\al$-forced by $n$ (i.e. $ \al^- [n^{\da}\2 ] \inc m^{\da} $); while
$(\al ,\be )$ is concept continuous iff it reflects incidence,
$n$ holds for $\al (g)$ whenever $g$ has each attribute $m$ that 
$\be$-generalizes $n$ (i.e. $n^{\da}\inc \be (m)^{\da}$), and dually,
$\be (m)$ holds for $h$ whenever $m$ holds for each  $\al$-specialization 
$g$ of $h$ (i.e. $h^{\ua} \inc \al (g)^{\ua}$).
 
\vx

A rather surprising consequence of the previous results is now:

\BT
A mapping pair $(\al ,\be )$ is a dense context embedding iff it is both 
conceptual and concept continuous.
\ET

\BP
Suppose $(\al ,\be )$ is a dense embedding. In order to show that it is
conceptual, use Lemmas \ref{dense}: for 
$h\nJ n$, choose $g,m$ such that 
$\al (g)\nJ \be (m)$, $h^{\ua}\inc \al (g)^{\ua}$
and $n^{\da}\inc \be (m)^{\da}$. 
Then $g \nI m$ since $(\al ,\be )$ preserves incidence, and
$\al^- [n^{\da}\2 ]\inc \al^- [\be (m)^{\da}\2 ]\inc m^{\da}$ since
$(\al ,\be )$ reflects incidence; dually, we get 
$\be^- [h^{\ua}\2 ]\inc g^{\ua}$, and Lemma \ref{conceptual} applies.
%
Now to concept continuity.  If $h\nJ \be (m)$, then by Lemma \ref{dense},
there is a $g$ with $h^{\ua}\inc \al (g)^{\ua}$ 
and $\al (g)\nJ \be (m)$, whence $g\nI m$ (because $(\al ,\be )$ 
preserves incidence).
This and a dual clue show that $(\al ,\be )$ is concept continuous,
on account of Lemma \ref{concont}.

\vh

Conversely, let $(\al ,\be )$ be conceptual and concept continuous.
By Lemmas \ref{conceptual} and \ref{concont}, $(\al ,\be )$ preserves and 
reflects incidence, so it is an embedding. 
For $h\nJ n$, there is an $m$ with $h\nJ \be (m)$ and 
$\al^- [n^{\da}\2 ]\inc m^{\da}$ (see the proof of Lemma \ref{conceptual}).
By Lemma \ref{concont}, we find a $g$ with $g\nI m$ and 
$h^{\ua}\inc\al (g)^{\ua}$. The assumption $\al (g) \2 J \2 n$ leads to the
contradiction $g\in \al^- [n^{\da}\2 ]\inc m^{\da}\3 $; hence $\al (g)\nJ n$. 
By Lemma \ref{dense}, we conclude that $\al$ (and similarly $\be$) is dense.
\EP  

\vh


\vx

Our final lemma shows that fullness or density of one partner in a
conceptual pair implies the corresponding property of the other.

\BL
{\rm (1)} If $(\al ,\be )$ is a context embedding then $\al$ and $\be$ are full.

\vh

\nid {\rm (2)} A conceptual pair $(\al ,\be )$ is an embedding 
iff $\al$ is full iff $\be$ is full.

\vh

\nid {\rm (3)} A conceptual pair $(\al ,\be )$ is dense iff $\al$ is dense
iff $\be$ is dense.
\EL

\BP
(1) If $(\al ,\be )$ is an embedding then $g \nI m$ implies 
$\al (g) \nJ \be (m)$, and 
$\al [m^{\da}\2 ]\inc \be (m)^{\da}$ (since incidence is preserved).
Hence, Lemma \ref{full} applies with $n=\be (m)$, showing that $\al$ is 
full, and by a dual argument, so is $\be$.

\vx

\nid (2) By (1) and Lemma \ref{conceptual}, it suffices to show that a
conceptual pair $(\al ,\be )$ for which $\al$ (or $\be$) is full must
reflect incidence. By Lemma \ref{full}, $g\nI m$ implies $\al (g)\nJ n$ 
for some $n\in \al [m^{\da}\2 ]^{\ua}$, and by Lemma \ref{conceptual}, 
this entails that $n$ belongs to $\be (m)^{\da\ua}$, i.e. 
$\be (m)^{\da}\inc n^{\da}$, and consequently $\al (g)\nJ \be (m)$.

\vx

\nid (3) Use Lemma \ref{insur} and the fact that the maps 
$\al^{\la}$ and $\be^{\la}$ agree up to the dual isomorphism between 
extent and intent lattices (see Section 6). 
\EP

\newpage

\section{Complete Lattices as Contexts}

We come now to the central part of our investigations, demonstrating that
our choice of morphisms was the ``right one'' from a categorical point of 
view. Each of the previously introduced classes of mapping pairs is closed 
under (componentwise) composition and may, therefore, serve as the morphism
class of a category of contexts. Specifically, we have 
the following categories of contexts and complete lattices
(see the next page for comments):

\vspace{2ex}

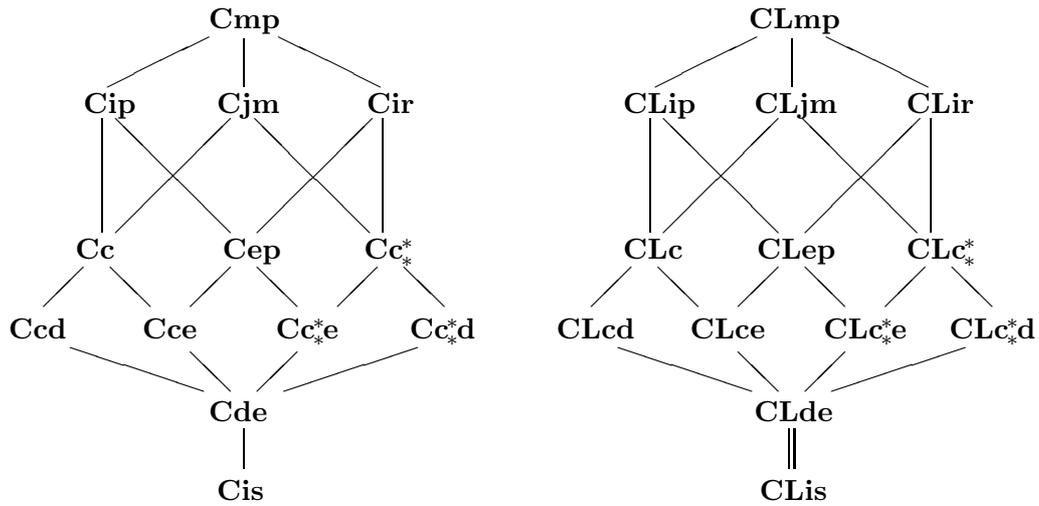
\nid {\bf Table 4.1}

\vspace{6ex}

\begin{picture}(0,420)
\put(-42,320){

\BTA {|l|l||l|l|}
\hline
{\em contexts$\!\!$}&{\em morphisms} &{\em lattices}$\!\!$&{\em morphisms}\\
\hline
\hline
{\bf Cmp}      &mapping pairs       &{\bf CLmp}$\!\!$&mapping pairs       \\
{\bf Cip}      &incidence preserving pairs&{\bf CLip}    &order preserving pairs\\
${\bf Cir}$  &incidence reflecting pairs&${\bf CLir}$&order reflecting pairs\\
{\bf Cep}      &context embedding pairs  &{\bf CLep}    &order embedding pairs\\
{\bf Cjm}&separately continuous pairs &{\bf CLjm}&join-meet preserving pairs$\!\!$\\
\hline
{\bf Cc}      &conceptual pairs    &{\bf CLc} &complete homomorphisms\\
${\bf Cc\!\2_*^*}$&concept continuous pairs&${\bf CLc\!\2_*}$&
doubly residuated maps\\
& &${\bf CLc\!\2^*}$&doubly residual maps\\
{\bf Ccd}     &dense conceptual pairs&${\bf CLcd}$ &
dense$\5$compl.$\5$homomorphisms$\!\!$\\
${\bf Cc\!\2_*^*\!\2 d}$ &dense$\5$concept$\5$continuous$\5$pairs$\!\!$&
${\bf CLc\!\2_*\!\2 d}\!\!$&doubly$\5$residuated$\5$surjections\\
& &${\bf CLc\!\2^*\!\2 d}\!\!$&doubly residual surjections\\
{\bf Cce}     &conceptual embeddings&{\bf CLce}  &complete embeddings\\
${\bf Cc\!\2_*^*\!\2 e}$ &concept$\5$continuous$\5$embeddings$\!\!$&
${\bf CLc\!\2_*\!\2 e}\!\!$&doubly$\5$residuated$\5$embeddings$\!\!$\\
& &${\bf CLc\!\2^*\!\2 e}\!\!$&doubly residual embeddings$\!\!$\\
{\bf Cde}     &dense embeddings &{\bf CLde} &dense$\6$complete$\6$embeddings\\
{\bf Cis}     &context isomorphisms &{\bf CLis}  &lattice isomorphisms\\
\hline
\ETA
}
\put(145,0){
\begin{picture}(200,170)

\put(101,-1){{\bf CLis}}
\put(112,11){\line(0,1){15}}
\put(114,11){\line(0,1){15}}

\put(99,29){{\bf CLde}}
\put(98,40){\line(-3,1){50}}
\put(108,40){\line(-1,1){15}}
\put(118,40){\line(1,1){15}}
\put(128,40){\line(3,1){50}}

\put(25,60){${\bf CLcd}$}
\put(38,72){\line(1,1){15}}

\put(75,60){{\bf CLce}}
\put(78,72){\line(-1,1){15}}
\put(93,72){\line(1,1){15}}

\put(125,60){${\bf CLc\!\2_*^*\!\2 e}$}
\put(148,72){\line(1,1){15}}
\put(133,72){\line(-1,1){15}}

\put(172,60){${\bf CLc\!\2_*^*\!\2 d}$}
\put(188,72){\line(-1,1){15}}

\put(50,90){{\bf CLc}}
\put(60,100){\line(0,1){43}}
\put(65,100){\line(1,1){43}}

\put(100,90){{\bf CLep}}
\put(108,100){\line(-1,1){43}}
\put(118,100){\line(1,1){43}}

\put(156,90){${\bf CLc\!\2_*^*}$}
\put(165,100){\line(0,1){43}}
\put(160,100){\line(-1,1){43}}

\put(50,145){{\bf CLip}}
\put(62,155){\line(2,1){38}}

\put(99,145){{\bf CLjm}}
\put(113,155){\line(0,1){18}}

\put(156,145){${\bf CLir}$}
\put(164,155){\line(-2,1){38}}

\put(97,177){{\bf CLmp}}

\end{picture}
}

\put(-60,0){
\begin{picture}(200,190)

\put(103,-1){{\bf Cis}}
\put(113,11){\line(0,1){15}}

\put(100,29){{\bf Cde}}
\put(98,40){\line(-3,1){50}}
\put(108,40){\line(-1,1){15}}
\put(118,40){\line(1,1){15}}
\put(128,40){\line(3,1){50}}

\put(25,60){{\bf Ccd}}
\put(38,72){\line(1,1){15}}

\put(75,60){{\bf Cce}}
\put(78,72){\line(-1,1){15}}
\put(93,72){\line(1,1){15}}

\put(125,60){${\bf Cc\!\2_*^*\!\2 e}$}
\put(148,72){\line(1,1){15}}
\put(133,72){\line(-1,1){15}}

\put(175,60){${\bf Cc\!\2_*^*\!\2 d}$}
\put(188,72){\line(-1,1){15}}

\put(50,90){{\bf Cc}}
\put(60,100){\line(0,1){43}}
\put(65,100){\line(1,1){43}}

\put(105,90){{\bf Cep}}
\put(108,100){\line(-1,1){43}}
\put(118,100){\line(1,1){43}}

\put(158,90){${\bf Cc\!\2_*^*}$}
\put(165,100){\line(0,1){43}}
\put(160,100){\line(-1,1){43}}

\put(53,145){{\bf Cip}}
\put(62,155){\line(2,1){38}}

\put(103,145){{\bf Cjm}}
\put(113,155){\line(0,1){18}}

\put(159,145){${\bf Cir}$}
\put(164,155){\line(-2,1){38}}

\put(100,177){{\bf Cmp}}

\end{picture}
}

\end{picture}

\newpage

\nid In view of the intended correspondences between contexts and complete 
lattices, we have included here a few less common types of morphisms:
an {\em order preserving mapping pair} $(\al ,\be )$ between posets is 
characterized by the implication $x\leq y \!\imp\! \al (x)\leq \be (y)$, 
while {\em order reflecting pairs} are characterized by the reverse 
implication, and {\em order embedding pairs} by the corresponding 
equivalence. 
By a {\em join-meet preserving pair}, we mean a mapping pair $(\al ,\be )$ 
between complete lattices such that $\al$ 
preserves arbitrary joins and $\be$ preserves arbitrary meets. 

\vx

Each category listed under ${\bf CLjm}$, the category of complete lattices
and join-meet preserving pairs,
may be embedded in that category by obvious identifications: complete 
homomorphisms $\vi$ are identified with pairs $(\vi ,\vi )$, 
doubly residuated maps $\psi$ with pairs $(\psi ,\psi^{**})$, and doubly 
residual maps $\psi$ with pairs $(\psi_{**},\psi )$. 
Thereby, the category ${\bf CLc}$ with complete homomorphisms, is 
identified with that subcategory of ${\bf CLjm}$ whose morphism pairs 
$(\al ,\be )$ have equal components $\al = \be$, whereas both ${\bf CLc}_*$,
the category with doubly residuated morphisms, and ${\bf CLc}^*$, the 
category with doubly residual morphisms,
are identified with that subcategory ${\bf CLc^*_*}$ of ${\bf CLjm}$ whose 
morphisms $(\al ,\be )$ satisfy the equation $\al^* = \be _*$.

\vx

\nid For any complete lattice ${\rm L} = (L,\leq )$, the {\em complete 
context} $\C\2 {\rm L} = (L,L,\leq )$ is the greatest context whose concept 
lattice is isomorphic to the original lattice ${\rm L}$. Thus, we have a 
functor 
$
\C : {\bf CLjm} \map {\bf Cjm},
$ 
sending ${\rm L}$ to $\C\2 {\rm L}$ and acting identically on mapping pairs.
Indeed, a map $\al$ between complete lattices preserves arbitrary joins iff
it is residuated, i.e. preimages of principal ideals are principal ideals,
and these are just the extents of the associated contexts; and dually, a
map preserves arbitrary meets iff it intent continuous. Thus, up to 
identification of complete lattices ${\rm L}$ with their complete
contexts $\C\2 {\rm L}$, the category {\bf CLjm} may be regarded as a full
subcategory of {\bf Cjm}.  

\vh

Note that by Lemma \ref{dense}, a map $\vi$ between complete lattices 
is join dense (resp. meet dense) iff it is extent (resp. intent) dense as a 
map between the associated contexts. Furthermore, recall that a
join- or meet-preserving dense map is already surjective; in particular,
dense complete embeddings are already isomorphisms. In all, we have:

\BT
\label{covarC}
Up to the aforementioned canonical identifications, $\C$ may be regarded
as a covariant full embedding functor from each of the categories of
complete lattices on the right hand list of Table {\rm 4.1} into the 
corresponding category of contexts on the left hand list. In particular, 
associating with any complete homomorphism $\vi$ the pair 
$(\vi ,\vi)$, one obtains full embeddings

\vx

of {\bf CLc} in {\bf Cc}, of {\bf CLcd} in {\bf Ccd}, of {\bf CLce} in 
{\bf Cce}, and of {\bf CLis} in {\bf Cde}.

\vx

\nid Similarly, associating with any doubly residuated map $\vi$ 
the pair $(\vi ,\vi^{**})$ and with any doubly residual map $\psi$ 
the pair $(\psi_{**} ,\psi)$, one obtains full embeddings 

\vx

of ${\bf CLc_*}$ and ${\bf CLc^*}$ in ${\bf Cc^*_*}$, of 
${\bf CLc_{*\!}d}$ and ${\bf CLc^{*\!}d}$ in ${\bf Cc_*^{*\!}d}$, etc. 


\ET

\BP
It remains to verify the following facts about a mapping pair $(\al ,\be )$
between complete contexts $\C\2{\rm K}$ and $\C\2{\rm L}$:
\BIT
\IT[(1)] $(\al ,\be )$ is conceptual iff $\al = \be$ is a complete 
homomorphism,
\mx 
\IT[(2)] $(\al ,\be )$ is dense and conceptual iff $\al =\be$ is a 
surjective complete homomorphism,
\mx
\IT[(3)] $(\al ,\be )$ is a conceptual embedding iff $\al =\be$ is an 
injective complete homomorphism,
\mx
\IT[(4)] $(\al ,\be )$ is a dense embedding iff $\al =\be$ is an 
isomorphism,\\[-2mm]
\mx
\IT[(${\rm 1^*}$)] $(\al ,\be )$ is concept continuous iff
$\al^{**}=\be$ (hence $\al =\be_{**}$ is doubly residuated and $\be$ is 
doubly residual),
\mx
\IT[(${\rm 2^*}$)] $(\al ,\be )$ is concept continuous and dense iff
$\al^{**}=\be$ is surjective\\
iff $\al =\be_{**}$ is surjective, 
\mx
\IT[(${\rm 3^*}$)] $(\al ,\be )$ is a concept continuous embedding iff
$\al^{**}=\be$ is injective\\
iff $\al = \be_{**}$ is injective.
\EIT
Concerning (1), note first that for a conceptual morphism $(\al ,\be )$
between complete contexts, both $\al$ and $\be$ are continuous, whence 
$\al$ preserves joins and $\be$ preserves meets. Furthermore, by Lemma
\ref{conceptual}, we have $\al (x)^{\ua\da} = \be [x^{\ua}\2]^{\da}$,
which means $\al (x) =\Inf \be[x^{\ua}\2 ] =\be (x)$, because $\be$ is
isotone. Conversely, every complete homomorphism $\vi$ yields a
conceptual morphism $(\vi ,\vi )$ (see \cite{EDM}). 

\vh

The equivalences (2) and (3) are now immediate consequences of the remarks 
before the theorem.
Concerning (4), recall the important fact that every dense context 
embedding $(\al ,\be )$ is conceptual (and concept continuous), whence
in the present situation $\al = \be$ is a join- and meet-dense complete
embedding and consequently an isomorphism.

\vx

\nid For ${\rm (1^*)}$, suppose first that $(\al ,\be )$ is concept 
continuous. Then the equation $\al^- [x^{\da}\2 ] =
\be^- [x^{\da\ua}\2 ]^{\da} = \be^- [x^{\ua}\2 ]^{\da}$ 
(see Lemma \ref{concont}) yields a map
\BM
\vi: {\rm L}\map {\rm K}, \ \vi (x) = \Sup \al^- [x^{\da}\2 ] = 
\Inf\be^- [x^{\ua}\2 ]
\EM
which is upper adjoint to $\al$ and lower adjoint to $\be$, whence
$\al^{**}=\vi^* = \be$ and $\al = \vi_* = \be_{**}$.
Conversely, for a doubly residual map $\psi : {\rm K}\map {\rm L}$,
the lower adjoint $\vi = \psi_*$ is a complete homomorphism, and the pair 
$(\psi_{**},\psi ) = (\vi_*,\vi^*)$ is concept continuous, since each 
concept of $\C\2 {\rm L}$ has the form $(x^{\da},x^{\ua})$, and consequently
$(\psi_{**}^- [x^{\da}\2 ],\psi ^-[x^{\ua}\2 ]) = 
(\vi (x)^{\da},\vi (x)^{\ua})$ is a concept of $\C\2 {\rm K}$.

\vx

The remaining statements are obtained as before, using the remarks on
density and embedding properties of join- or meet-preserving maps.
\EP

Note that the above claims remain valid if we replace join- and meet-preserving 
maps with residuated and residual maps between arbitrary posets. Thus the category
of posets with residuated and residual maps is a full subcategory of the 
category ${\bf Cc}$, etc.

\newpage

\section{The Concept Lattice as a Covariant Functor}

The point is now that the embedding functor $\C$ has a left adjoint,
sending on the object level each context to its concept lattice. Thus,
for any context $\KK$, the concept lattice $\B\2  \KK$ may be viewed as the
{\em free complete lattice over} $\KK$.
On the morphism level, we define for any mapping pair $(\al ,\be )$ between
contexts $\KK = (G,M,I)$ and $\LL =(H,N,J)$ a ``lifted'' mapping pair 
$(\al^{\ra} ,\be^{\ra})$ by

\vx

\ \ \ $\al^{\ra}:\B\2 \KK\map \B\2 \LL ,\ \1\al^{\ra}(A,B)=(\al [A]^{\ua\da},
\al [A]^{\ua}),$

\vh

\ \ \ $\be^{\ra}:\B\2 \KK\map \B\2 \LL ,\ \be^{\ra}(A,B)=(\be [B]^{\da},
\be [B]^{\da\ua}).
\EM

\BPR
\label{lift}
A mapping pair $(\al ,\be )$ between contexts $\KK$ and $\LL$ is

\vh

- separately continuous iff $\al^{\ra}$ is join and $\be^{\ra}$ is 
meet preserving,


- conceptual iff $\al^{\ra}\! =\be^{\ra}$ is a complete homomorphism,


- dense conceptual iff $\al^{\ra}\! =\be^{\ra}$ is a surjective complete
homomorphism,


- a conceptual embedding iff $\al^{\ra}\! =\be^{\ra}$ is a complete embedding,


- a dense embedding iff $\al^{\ra}\! =\be^{\ra}$ is an isomorphism.

\vx

\nid On the other hand, $(\al ,\be )$ is concept continuous 

\vh

\ \ iff $\ \al^{\ra}$ is doubly residuated with $\al^{\ra **}\! =\be^{\ra}$

\ \ iff $\ \be^{\ra}$ is doubly residual with $\al^{\ra}\! = {\be^{\ra}}_{**}$.

\EPR

\BP
The first equivalence follows from Theorem \ref{cont}. The
second is obtained from the first one and Lemma \ref{conceptual} (cf. 
\cite{EDM}). For the statements about density and embedding properties, 
apply Lemma \ref{insur}.

\vx

In the case of a concept continuous pair $(\al ,\be )$, we have a 
well-defined complete homomorphism
\BM
\vi : \B\2 \LL \map \B\2 \KK ,\ (C,D)\mapsto (\al^- [C], \be^- [D])
\EM
because $\al^-$ and $\be^-$ preserve arbitrary intersections. 
The equivalences
\BM
\al^{\ra}(A,B)\leq (C,D) \!\iff\! 
\al [A]\inc C \!\iff\! A\inc \al^- [C] \!\iff\! (A,B)\leq \vi (C,D)
\EM
show that $\vi$ is the upper adjoint of $\al^{\ra}$. A dual argument shows
that $\vi$ is the lower adjoint of $\be^{\ra}$, whence 
$\al^{\ra **}= \vi^* = \be^{\ra}$ and $\al^{\ra} = \vi_* ={\be^{\ra}}_{**}$.
\\
Conversely, if $\al^{\ra **} =\be^{\ra}$ then $\vi = \al^{\ra *}$ is a 
complete homomorphism from $\C\2{\rm L}$ into $\C\2{\rm K}$ which is upper
adjoint to $\al^{\ra}$ and lower adjoint to $\be^{\ra}$. For any two 
concepts $(A,B)\in \B\2 \KK$ and $(C,D)\in \B\2 \LL$, we have the equivalences
\BM
(A,B)\leq \vi (C,D)\!\iff\! \al^{\ra} (A,B)\leq (C,D)
\!\iff\! (A,B)\leq (\al^- [C],\al^- [C]^{\ua}),$

\vh

\ \ \ $
\vi (C,D)\leq (A,B) \!\iff\! (C,D)\leq \be^{\ra}(A,B) 
\!\iff\! (\be^- [D]^{\da}, \be^- [D])\leq (A,B),
\EM
which amount to the equations
\BM
\vi (C,D) = (\al^- [C], \be^- [D]) = (\al^- [C],\al^- [C]^{\ua})= 
(\be^- [D]^{\da},\be^- [D]),
\EM
establishing the concept continuity of $(\al ,\be )$.
\EP

\vx

\nid Of basic importance for Formal Concept Analysis are the dense 
embeddings

\vx

\ \ \ $\eta_{\KK} = (\ga_{\KK} ,\mu_{\KK}) : \KK \map \C\B\2 \KK = 
(\B\2 \KK ,\B\2 \KK ,\leq )$, \ where 

\vx

\ \ \ $\ga_{\KK} = \ga : G\map \B\2  \KK\ $ and 
$\ \mu_{\KK}=\mu : M\map \B\2 \KK
\EM 
are the natural object and attribute embeddings, respectively (cf. 
\cite{W1}):
\BM
\ga (g) =  (g^{\ua\da},g^{\ua}) \ $ and $\ \mu (m) = (m^{\da}, m^{\da\ua}).
\EM
It was pointed out in the preprint version of \cite{EDM} that these 
embeddings are the reflection 
morphisms for a reflector $\C\B$ from the category {\bf Cc} of contexts and 
conceptual mapping pairs to the full subcategory {\bf CCc} of complete
contexts and complete homomorphisms. The construction from \cite{EDM} may be 
extended, in a straightforward manner, to the categories {\bf Cjm} and
{\bf CLjm$\3$}: for any separately continuous mapping pair $(\al ,\be )$ between
contexts $\KK$ and $\LL$, we put $\B (\al ,\be )=(\al^{\ra},\be^{\ra})$. 
In case $(\al ,\be )$ is conceptual, we have $\al^{\ra} =\be^{\ra}$, so that
we may replace $(\al^{\ra},\be^{\ra})$ with $\B (\al, \be) = \al^{\ra}$.
Similarly, if $(\al ,\be )$ is concept continuous, the second component
is determined by the first one via $\al^{\ra **}\! =\be^{\ra}$, so that 
it is more convenient to put $\B (\al, \be) = \al^{\ra}$ (or, alter\-nately,
$\B (\al, \be ) = \be^{\ra}$). Concerning adjoint functors, we refer to \cite{AHS}.

\BT
\label{covarB}
$\B$ is a covariant functor 

\vx

\BTA {lllllllll}
from &$\!\!{\bf Cjm}\!\!$&$\!\!$to$\!\!$&$\!\!{\bf CLjm}\!\!$ &
\hspace{3ex} &
&&&\\ 
from &$\!\!{\bf Cc}\!\!$&$\!\!$to$\!\!$&$\!\!{\bf CLc}\!\!$ &
\hspace{3ex} &
from &$\!\!{\bf Cc_*^*}\!\!$&$\!\!$to$\!\!$&$\!\!{\bf CLc_*}\6\ $ 
and $\ {\bf CLc^*}$\\
from &$\!\!{\bf Ccd}\!\!$&$\!\!$to$\!\!$&$\!\!{\bf CLcd}\!\!$ &
\hspace{3ex} &
from &$\!\!{\bf Cc_*^*d}\!\!$&$\!\!$to$\!\!$&$\!\!{\bf CLc_*d}$ 
and $\ {\bf CLc^*e}$\\
from &$\!\!{\bf Ccd}\!\!$&$\!\!$to$\!\!$&$\!\!{\bf CLce}\!\!$ &
\hspace{3ex} &
from &$\!\!{\bf Cc_*^*e}\!\!$&$\!\!$to$\!\!$&$\!\!{\bf CLc_*e}$
and $\ {\bf CLc^*e}$\\
from &$\!\!{\bf Cde}\!\!$&$\!\!$to$\!\!$&$\!\!{\bf CLis}\!\!$ &
\hspace{3ex} &
&&&

\ETA

\vh

\nid Furthermore, $\B$ is left adjoint to the functor $\C$ in
the opposite direction. The unit of the adjunctions is $\eta$, and the 
counit is an isomorphism.
\ET

\BP
The functor properties of $\B$ are easily checked, using Theorem \ref{cont}
and Proposition \ref{lift}. Concerning adjointness, we have to show that 
for every complete lattice ${\rm L}$ and for every separately continuous 
morphism $(\al ,\be )$ from an arbitrary context $\KK = (G,M,I)$ 
into $\C\2{\rm L} =(L,L,\leq)$, there is a unique pair 
$(\al^{\ve},\be^{\we})$ such that 

\vh

\ \ \ $\al^{\ve} : \B\2 \KK \map {\rm L}$ preserves joins,
$\be^{\we} : \B\2 \KK \map {\rm L}$ preserves meets, and 
\BM
(\al ,\be ) = (\al^{\ve},\be^{\we})\circ \eta_{\KK}\3$, i.e. $
\al =\al^{\ve}\!\circ\ga_{\KK}\3$ and $\be =\be^{\we}\!\circ\mu_{\KK}.
\EM
Defining for each concept $(A,B)\in\B\2 \KK$
\BM
\al^{\ve} (A,B)=\Sup \al [A]$ \ and \ $\be^{\we}(A,B)= \Inf \be [B],
\EM
we see that $\al ^{\ve}$ is lower adjoint to the map
\BM
\al^{\ve *} : {\rm L} \map \B\2 \KK , \ x \mapsto 
(\al^-[{\da x}],\al^- [{\da x}]^{\ua}),
\EM 
and that 
$\be ^{\we}$ is upper adjoint to the map
\BM
{\be^{\we}}_* : {\rm L} \map \B\2 \KK , \ x \mapsto 
(\be^-[{\ua x}]^{\da},\be^- [{\ua x}]).
\EM 
Hence, $\al^{\ve}$ preserves joins, $\be^{\we}$ preserves meets, and
\BM
\al^{\ve}\!\circ \ga (g) \ = \Sup \al [\3 g^{\ua\2\da}\2 ]\3 = \3\al (g)\6$ 
by continuity of $\al$,

\vh

\ \ \ $
\be^{\we}\!\circ \mu (m) = \Inf \be [m^{\da\ua}\2 ] = \be (m)$ by continuity
of $\be .
\EM
The uniqueness of a join-meet preserving pair $(\al^{\ve},\be^{\we})$
with $\al =\al^{\ve}\!\circ\ga_{\KK}\3$ and $\3\be =\be^{\we}\!\circ\mu_{\KK}$
follows from join-density of $\ga_{\KK}$ and meet-density of $\mu_{\KK}$.

\vh

Now, if $(\al ,\be ) : \KK \map \LL$ is an arbitrary {\bf Cjm}-morphism 
(that is, a separately continuous mapping pair) then there is a unique 
{\bf CLjm}-morphism 
$(\vi ,\psi ) : \B\2 \KK \map \B\2 \LL$ satisfying the identity 
$
(\vi ,\psi )\circ\eta_{\KK} = \eta_{\LL}\circ (\al ,\be ),
$
namely
\BM
\vi (A,B) = \Sup \ga_{\LL}[\al [A\2 ]] = 
(\al [A\2 ]^{\ua\da},\al [A\2 ]^{\ua}) =
\al^{\ra} (A,B),$

\vx

\ \ \ $
\psi (A,B) = \Inf \mu_{\LL}[\be [B]] = (\be [B]^{\da\ua},\be [B]^{\da}) =
\be^{\ra} (A,B),
\EM
and consequently $(\vi ,\psi ) = \B (\al ,\be )$.
This shows that $\B$ is in fact left adjoint to $\C$, and that $\eta$ is
the unit of the adjunction. The counit $\ep$ is constituted by 
the natural isomorphisms

\vx

\ \ \ $
\ep_{{\rm L}} :\B\C\2 {\rm L} \map {\rm L}$ \ with \ 
$\ep_{\rm L}({\ua x},{\da x}) = x.$
\EP

\BPI (300,100)
\put(0,80){$\KK$}
\put(13,84){\vector(1,0){60}}
\put(35,88){$\al$}
\put(40,73){$\be$}
\put(77,80){$\C\2{\rm L}$}

\put(4,77){\vector(0,-1){60}}
\put(-11,49){$\ga_{\KK}$}
\put(7,43){$\mu_{\KK}$}

\put(-8,5){$\C\B\2 \KK$}
\put(12,18){\vector(1,1){60}}
\put(32,46){$\al^{\ve}$}
\put(43,40){$\be^{\we}$}

\put(127,80){${\rm L}$}
\put(41,5){$\B\2 \KK$}
\put(120,76){\vector(-1,-1){60}}
\put(76,46){$\al^{\ve *}$}
\put(91,40){${\be^{\we\!}}_*$}

\put(120,10){\vector(-1,0){60}}
\put(87,13){$\al^-$}
\put(93,-2){$\be^-$}
\put(133,43){$\ep_{\rm L}$}

\put(122,5){$\B\C\2{\rm L}$}
\put(130,16){\vector(0,1){60}}

\put(200,80){$\KK$}
\put(213,84){\vector(1,0){60}}
\put(237,88){$\al$}
\put(243,73){$\be$}
\put(277,80){$\LL$}

\put(204,77){\vector(0,-1){60}}
\put(189,49){$\ga_{\KK}$}
\put(208,43){$\mu_{\KK}$}
\put(186,5){$\C\B\2 \KK$}

\put(279,77){\vector(0,-1){60}}
\put(264,49){$\ga_{\LL}$}
\put(282,43){$\mu_{\LL}$}
\put(274,5){$\C\B\2 \LL$}
\put(213,10){\vector(1,0){60}}
\put(237,14){$\vi$}
\put(243,0){$\psi$}

\EPI

\vx

\nid Invoking Theorem \ref{cont} once more, we conclude:

\BCO
\label{char}
A mapping pair $(\al ,\be ) : \KK \map \LL$ is separately continuous iff
there exists a (unique) {\bf CLjm}-morphism 
$(\vi ,\psi ) : \C\B\2 \KK \map \C\B\2 \LL$
with
\BM
\vi\circ \ga_{\KK} = \ga_{\LL}\circ\al$ and 
$\psi\circ \mu_{\KK} = \mu_{\LL}\circ\be \2 ,
\EM
namely $(\vi ,\psi ) = (\al^{\ra},\be^{\ra})$. Furthermore,

\vh

$(\al ,\be )$ is conceptual iff $\vi =\psi$ is a complete homomorphism,

$(\al ,\be )$ is concept continuous iff $\vi^*=\psi_*$ is a complete homomorphism,

$(\al ,\be )$ is a dense embedding iff $\vi = \psi \ (=\psi_{**})$ is an
isomorphism.
\ECO

\BCO
Up to identification between complete lattices and complete contexts,
we have full reflective subcategories

\vx

\BTA {llll}

$\hspace{-1.2ex}{\bf CLjm}\hookrightarrow {\bf Cjm},$&${\!\bf CLc}\2\hookrightarrow {\bf Cc},$
&$\!{\bf CLcd}\2\hookrightarrow {\bf Ccd},$&$\!{\bf CLce}\2\hookrightarrow {\bf Cce}, 
$\\[2mm]
$\hspace{-1.2ex}{\bf CLis}\ \hookrightarrow \4{\bf Cde},$&
$\!{\bf CLc\!\2_*}\!\!\hookrightarrow {\bf Cc\!\2_*^*},$&
$\!{\bf CLc\!\2_*\!\2 d}\!\hookrightarrow {\bf Cc\!\2 _*^*\!\2 d},$&
$\!{\bf CLc\!\2_*\!\2 e}\!\hookrightarrow {\bf Cc\!\2 _*^*\!\2 e},$
\ETA 

\vx

\nid and analogous reflections for the corresponding categories with
doubly residual morphisms.
\ECO

\BCO
\label{charc}
For all contexts $\KK$ and complete lattices ${\rm L}$, a mapping pair
$(\al ,\be ): \KK \map\C\2{\rm L}$ is conceptual iff there is a unique 
complete homomorphism \\
$\vi : \B\2 \KK \map {\rm L}$ such that
$(\al ,\be ) = (\vi ,\vi )\circ\eta_{\KK}$, \ i.e.\ 
$\al =\vi\circ\ga_{\KK}$ and $\be = \vi\circ\mu_{\KK}$.
Moreover, $(\al ,\be )$ is a dense embedding iff $\vi$ is an isomorphism.
\ECO

\BPI (300,100)
\put(50,80){$\KK$}
\put(63,84){\vector(1,0){60}}
\put(80,88){$(\al ,\be )$}
\put(127,80){$\C\2{\rm L}$}

\put(54,77){\vector(0,-1){60}}
\put(-18,46){$\eta_{\KK}=(\ga_{\KK} ,\mu_{\KK} )$}

\put(42,5){$\C\B\2 \KK$}
\put(62,18){\vector(1,1){60}}
\put(100,46){$(\vi ,\vi )$}

\put(177,80){${\rm L}$}
\put(91,5){$\B\2 \KK$}
\put(110,16){\vector(1,1){60}}
\put(148,46){$\vi$}

\put(227,80){${\rm L}$}
\put(141,5){$\B\2 \KK$}
\put(220,76){\vector(-1,-1){60}}
\put(180,46){$\vi^*$}
\put(191,40){$\vi_*$}

\put(220,10){\vector(-1,0){60}}
\put(187,13){$\al^{\ra *}$}
\put(193,-2){${\be^{\ra}\!}_*$}
\put(233,43){$\ep_{\rm L}$}

\put(222,5){$\B\C\2{\rm L}$}
\put(230,16){\vector(0,1){60}}

\EPI

\nid These results extend known facts about the Dedekind-MacNeille 
completion (by cuts) of ordered sets (cf. \cite{EOV,EDM}). Another 
immediate application is the {\em Fundamental Theorem on Concept Lattices}
(see e.g. \cite{W1}), saying that a concept lattice 
$\B (G,M,I)$ is isomorphic to a given complete lattice ${\rm L}$ iff there 
exists a join-dense map $\ga$ from $G$ into ${\rm L}$ and a meet-dense map
$\mu$ from $M$ into ${\rm L}$ such that $g\2 I\2 m \iff\ga (g)\leq\mu (m)$. 

\vspace{1ex}

If $(G,M,I)$ is a subcontext of a context $(H,N,J)$ 
(that is, $G\!\inc\! H$, $M\!\inc\! N$ and $I = (G\!\times\! M)\cap J$), 
we may consider the respective inclusion maps.

\BCO
\label{subc}
For a subcontext $\KK = (G,M,I)$ of a context $\LL = (H,N,J)$,
the following conditions are equivalent:
\BIT
\aa The inclusion maps $G\hookrightarrow H$ and $M\hookrightarrow N$
form a conceptual pair.
\bb $({A^J\!}\cap M)_J \inc {A^J\!}_J$ for $A\inc G\ $ and 
$\ ({B_J\!}\cap G)^J \inc {B_J}^J$ for $B\inc M$.
\cc If $h\nJ n$ then there are $g\nI m$ with $h^J\cap M\inc g^J$ and 
$n_J\cap G\inc m_J$.
\dd $(A,B)\mapsto ({A^J\!\!}_J,A^J\!\2 ) = (B_J,{B_{\!\2 J}\!\3}^J\!\2)$ 
is a complete homomorphism from $\B\KK$ to $\B\LL$.
\ee There is a complete homomorphism $\vi :\B\2 \KK \map \B\2 \LL$ with
$\vi\circ\eta_{\KK} = \eta_{\LL}$.
\EIT
\ECO

Similarly, taking for $\al$ and $\be$ identity maps but different incidence
relations, we arrive at

\BCO
\label{rel}
For two contexts $\KK = (G,M,I)$ and $\LL = (G,M,J)$ with the same 
underlying sets, the following conditions are equivalent:
\BIT
\aa The identity pair $(id_G ,id_M )$ is conceptual.
\bb ${A^J\!}_J = {A^I\!}_J$ for $A\inc G\ $ and 
$\ {B_J}^J = {B_I}^J$ for $B\inc M$.
\cc $I\inc J$, and for all $h\nJ n$ there are $g\nI m$ with 
$h^J\inc g^I$ and $n_J\inc m_I$.
\dd , \ {\rm (e)} \ As in Corollary {\rm \ref{subc}}.
\EIT 
\ECO

Note that in contrast to Corollary \ref{subc}, the complete homomorphisms in
Corollary \ref{rel} are always surjective, because identity maps are trivially
dense (but not necessarily full).

\section{The Concept Lattice as a Contravariant Functor}

We have already seen in the preceding section that every concept continuous
context morphism $(\al ,\be ) : \KK \map \LL$ gives rise to a complete
homomorphism in the opposite direction,
\BM
\vi : \B\2 \LL \map \B\2 \KK, \ (C,D)\mapsto (\al^- [C],\be^- [D]).
\EM
More generally, any separately continuous context morphism 
$(\al ,\be ) : \KK\map \LL$ induces a meet-preserving map
\BM
\al^{\la} = \al^{\ra *} :\B\2 \LL \map \B\2 \KK ,\ (C,D)\mapsto (\al^- [C],\al^- [C]^{\ua})
\EM
and a join-preserving map
\BM
\be^{\la} = {\be^{\ra}\!}_*:\B\2 \LL \map \B\2 \KK ,\ (C,D)\mapsto (\be^- [D]^{\da},\be^- [D]),
\EM
and these two maps coincide iff $(\al ,\be )$ is concept continuous.

\vx

The category {\bf CLmj} of complete lattices and meet-join preserving pairs
$(\psi ,\vi )$ (where $\psi : {\rm K} \map {\rm L}$ preserves arbitrary meets
and $\vi : {\rm K} \map {\rm L}$ preserves arbitrary joins) is isomorphic
to the category {\bf CLjm} by means of two essentially different functors:
one of them exchanges the first and the second component in the mapping
pairs, while the other keeps the morphisms fixed and reverses the lattice
orders. But there is also a {\em dual} isomorphism between {\bf CLjm} and
{\bf CLmj}, obtained by passing to the order-theoretical adjoints:

\vx

\BTA {lll}
$\4\G : {\bf CLjm} \map {\bf CLmj},$&$\4\G {\rm L} ={\rm L},$&
$\4\G (\al, \be ) = (\al^*,\be_*),$\\
$\H : {\bf CLmj} \map {\bf CLjm},$&$\H {\rm L} ={\rm L},$&
$\H (\psi, \vi ) = (\psi_*,\vi^*).$
\ETA

\vx

\nid Obviously, these mutually inverse functors induce dual
isomorphisms between the categories ${\bf CLc}$ and ${\bf CLc^*_*}$ etc.

\vx

Now, by our previous considerations, we have a contravariant functor
\BM
\B^{\la} : {\bf Cjm} \map {\bf CLmj}\2 , \ \B^{\la}\KK = \B\2 \KK\2 ,
\ \B^{\la} (\al ,\be ) = (\al^{\la} ,\be^{\la} ),
\EM
and also a contravariant functor in the other direction, namely

\BM
\C^{\la} : {\bf CLmj} \map {\bf Cjm}\2 ,\5\ \C^{\la}{\rm L} =\C\2{\rm L}\2 ,
\ \ \C^{\la} (\psi ,\vi ) = (\psi_*, \vi^*).
\EM
Moreover, these functors are linked with the covariant functors $\B$ and $\C$
by the identities
\BM
\B^{\la} = \G\circ\B, \ \B = \H\circ\B^{\la}, \ \C^{\la} =\C\circ\H,
\ \C = \C^{\la}\!\circ\G .
\EM 
Therefore, the adjunction in Theorem \ref{covarB} turns into a dual 
adjunction for the corresponding contravariant functors:

\BT
\label{contraB}
The contravariant functor $\B^{\la} : {\bf Cjm}\map {\bf CLmj}$ is dually
adjoint to the contravariant functor $\C^{\la} : {\bf CLmj}\map {\bf Cjm}$.
Furthermore, these functors induce dual adjunctions between
\BM
{\bf Cc}$ and ${\bf CLc\!\2_*^*}$, \ ${\bf Cc\!\2_*^*}$ and ${\bf CLc}$, \
${\bf Cde}$ and ${\bf CLis}$, etc.
\ET

\BCO
\label{charcon}
For all contexts $\KK$ and complete lattices ${\rm L}$, a mapping pair 
$(\al ,\be ) :\KK \map \C\2{\rm L}$ is concept continuous iff there is a
unique complete homomorphism $\vi : {\rm L} \map \B\2 \KK$ such that
$\vi_*\circ\ga_{\KK} = \al$ and $\vi^*\circ\mu_{\KK} = \be$.
\ECO

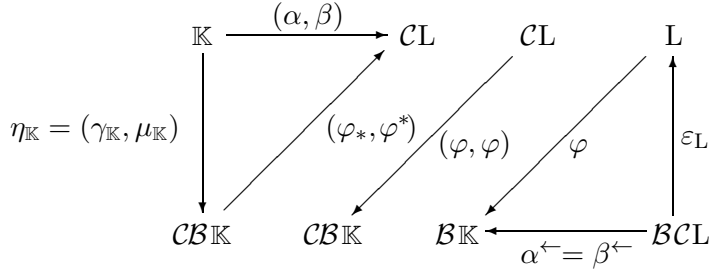
\BPI (300,100)
\put(50,80){$\KK$}
\put(63,84){\vector(1,0){60}}
\put(80,88){$(\al ,\be )$}
\put(127,80){$\C\2{\rm L}$}

\put(54,77){\vector(0,-1){60}}
\put(-18,46){$\eta_{\KK}=(\ga_{\KK} ,\mu_{\KK} )$}

\put(42,5){$\C\B\2 \KK$}
\put(62,18){\vector(1,1){60}}
\put(100,46){$(\vi_* ,\!\2\vi^*\! )$}

\put(172,80){$\C\2{\rm L}$}
\put(91,5){$\C\B\2 \KK$}
\put(170,76){\vector(-1,-1){60}}
\put(142,40){$(\vi ,\vi )$}

\put(227,80){${\rm L}$}
\put(141,5){$\B\2 \KK$}
\put(220,76){\vector(-1,-1){60}}
\put(191,40){$\vi$}

\put(220,10){\vector(-1,0){60}}
\put(173,-2){$\al^{\la}\!\! =\be^{\la}$}
\put(233,43){$\ep_{\rm L}$}

\put(222,5){$\B\C\2{\rm L}$}
\put(230,16){\vector(0,1){60}}

\EPI

\BCO
A mapping pair $(\al ,\be )$ between contexts $\KK$ and $\LL$ is concept 
continuous iff there exists a (unique) complete homomorphism 
$\vi : \B\2 \LL \map \B\2 \KK$
such that $\vi_*\circ\ga_{\KK} = \ga_{\LL}\circ\al$ and 
$\vi^*\circ\mu_{\KK} = \mu_{\LL}\circ\be\2$, namely $\vi = \al^{\la}\! =\be^{\la}$.
\ECO




\BCO
\label{subcontext}
For a subcontext $\KK = (G,M,I)$ of a context $\LL = (H,N,J)$, the following 
conditions are equivalent:
\BIT
\aa The inclusion morphism from $\KK$ into $\LL$ is concept continuous.
\bb $\KK$ is compatible, that is,
$(C^J\cap M)_J\cap G\inc {C^J}\!\!_J$ for all $C\!\2\inc\!\2 H$ and\\
$(D_{\!\2 J}\cap G)^J\cap M\inc {D_{\!\2 J}}^J$ for all $D\!\2\inc\!\2 N\2$.
\cc For all $g\!\in\! G$ and $n\!\in\! N\!\s- g^J$ there is an 
$m\!\in\! M\!\s- g^J$ with $n_J\inc m_J$, and\\
for all $m\!\2\in\!\2 M\!$ and $h\!\in\! H\!\s- m_J$ there is a
$g\!\in\! G\!\s- m_J$ with $h^J\inc g^J$.
\dd The trace map $\vi : \B\2 \LL\map \B\2 \KK, \ (C,D)\mapsto (C\cap G,D\cap M)$
is a complete homomorphism.
\ee There exists a unique complete homomorphism $\vi$ from $\B\2 \LL$ onto 
$\B\2 \KK$ with $\vi_*\circ\ga_{\KK\!} = \ga_{\LL}$ and 
$\vi^*\circ\mu_{\KK\!}=\mu_{\LL}$.
\EIT
\ECO



Another frequently used application of our results on concept continuous maps is 
obtained by taking for $\al$ and $\be$ identity maps.

\BCO
\label{closed}
For contexts $\KK = (G,M,I)$ and $\LL = (G,M,J)$, the following
conditions are equivalent:
\BIT
\aa The identity pair $(id_G,id_M): \KK\map \LL$ is concept continuous.
\bb $J\inc I, \ {A^J\!}_I \inc {A^J\!}_J$ for $A\inc G$, and ${B_J}^I\inc {B_J}^J$
for $B\inc M$.
\cc $J\inc I$, and $(g,m)\in I\s- J$ implies $h\nI m$ for some $h\in G$ with 
$g^J\inc h^J$ and $g\nI n$ for some $n\in M$ with $m_J\inc n_J$.
\dd $J$ is a closed relation of $\KK$, that is, 
$J\inc G\!\times\! M$ and $\B\2 \LL\inc \B\2\KK$. 
\ee $\B\2  \LL$ is a complete sublattice of $\B\2  \KK$.
\EIT
\ECO

A great part of these two corollaries has been discovered earlier by Ganter, Wille
and Reuter (see \cite{W2,W3} and \cite{RW}).


\newpage


\vspace{2ex}
\BC
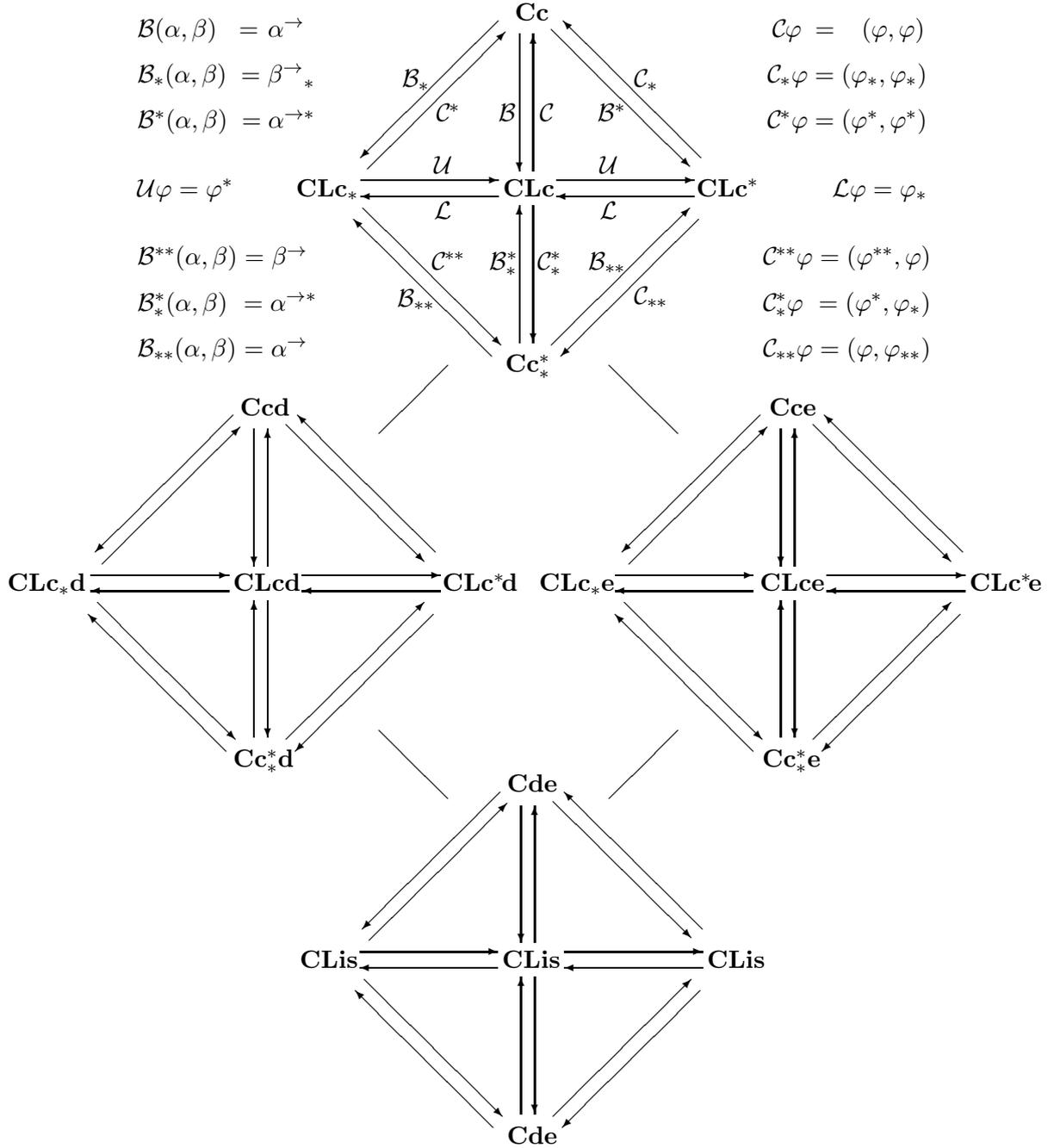
{\bf \ Diagram. \ Categories of Contexts and Complete Lattices}
\EC

\vx

\BPI (350,170)(-80,0)

\put(91,81){${\bf CLc}$}
\put(1,81){${\bf CLc_*}$}
\put(175,81){${\bf CLc^*}$}
\put(96,158){${\bf Cc}$}
\put(92,4){${\bf Cc^*_*}$}

\put(89,114){$\B$}
\put(106,114){$\C$}

\put(132,114){$\B^*$}
\put(61,114){$\C^{*}$}

\put(47,129){$\B_*$}
\put(147,129){$\C_{*}$}

\put(85,49){$\B_*^*$}
\put(105,49){$\C_*^*$}

\put(147,34){$\C_{**}$}
\put(128,49){$\B_{**}$}

\put(59,49){$\C^{**}$}
\put(45,33){$\B_{**}$}

\put(61,91){$\U$}
\put(61,71){$\L$}
\put(133,91){$\U$}
\put(133,71){$\L$}

\put(-68,151){$\B (\al ,\be )\ \6 =\2  \al^{\ra}$}
\put(-68,131){$\B_*(\al ,\be )\6 = {\be^{\ra}}_*$}
\put(-68,111){$\B^*(\al ,\be )\6 = \al^{\ra *}$}

\put(-68,81){$\U \vi = \vi^*$}

\put(-68,51){$\B^{**} (\al ,\be )= \be^{\ra}$}
\put(-68,31){$\B^*_*(\al ,\be )\ = \al^{\ra *}$}
\put(-68,11){$\B_{**} (\al ,\be )= \al^{\ra}$}

\put(207,151){$\C \vi \6 = \ \ \3 (\vi ,\vi )$}
\put(205,131){$\C_{*}\vi = (\vi_* ,\vi_*)$}
\put(205,111){$\C^{*}\vi = (\vi^* ,\vi^*)$}

\put(234,81){$\L \vi= \vi_*$}

\put(203,51){$\C^{**}\vi = (\vi^{**},\vi )$}
\put(203,31){$\C_*^*\vi \ = (\vi^*,\vi_*)$}
\put(203,11){$\C_{**}\vi = (\vi ,\vi_{**})$}

\put(114,88){\vector(1,0){60}}
\put(174,81){\vector(-1,0){60}}
\put(29,88){\vector(1,0){60}}
\put(89,81){\vector(-1,0){60}}
\put(98,17){\vector(0,1){60}}
\put(104,77){\vector(0,-1){60}}
\put(104,92){\vector(0,1){60}}
\put(98,152){\vector(0,-1){60}}

\put(112,16){\vector(1,1){60}}
\put(176,71){\vector(-1,-1){60}}

\put(177,98){\vector(-1,1){60}}
\put(112,154){\vector(1,-1){60}}

\put(32,93){\vector(1,1){60}}
\put(89,158){\vector(-1,-1){60}}

\put(86,11){\vector(-1,1){60}}
\put(30,76){\vector(1,-1){60}}

\EPI

\begin{picture}(300,300)

\put(105,290){\line(1,1){30}}

\put(105,160){\line(1,-1){30}}

\put(235,290){\line(-1,1){30}}

\put(235,160){\line(-1,-1){30}}

\put(-51,140)
{

\begin{picture}(100,100)

\put(89,81){${\bf CLcd}$}
\put(-9,81){${\bf CLc_*d}$}
\put(180,81){${\bf CLc^{*\!}d}$}
\put(92,158){${\bf Ccd}$}
\put(89,4){${\bf Cc^{*\!}_*d}$}

\put(119,88){\vector(1,0){60}}
\put(179,81){\vector(-1,0){60}}
\put(27,88){\vector(1,0){60}}
\put(87,81){\vector(-1,0){60}}
\put(98,17){\vector(0,1){60}}
\put(104,77){\vector(0,-1){60}}
\put(104,92){\vector(0,1){60}}
\put(98,152){\vector(0,-1){60}}

\put(112,16){\vector(1,1){60}}
\put(176,71){\vector(-1,-1){60}}

\put(177,98){\vector(-1,1){60}}
\put(112,154){\vector(1,-1){60}}

\put(32,93){\vector(1,1){60}}
\put(89,158){\vector(-1,-1){60}}

\put(86,11){\vector(-1,1){60}}
\put(30,76){\vector(1,-1){60}}

\end{picture}
}

\put(178,140)
{

\begin{picture}(100,100)

\put(89,81){${\bf CLce}$}
\put(-7,81){${\bf CLc_*e}$}
\put(180,81){${\bf CLc^{*\!}e}$}
\put(93,158){${\bf Cce}$}
\put(90,4){${\bf Cc^{*\!}_*e}$}

\put(118,88){\vector(1,0){60}}
\put(178,81){\vector(-1,0){60}}
\put(26,88){\vector(1,0){60}}
\put(86,81){\vector(-1,0){60}}
\put(98,17){\vector(0,1){60}}
\put(104,77){\vector(0,-1){60}}
\put(104,92){\vector(0,1){60}}
\put(98,152){\vector(0,-1){60}}

\put(112,16){\vector(1,1){60}}
\put(176,71){\vector(-1,-1){60}}

\put(177,98){\vector(-1,1){60}}
\put(112,154){\vector(1,-1){60}}

\put(32,93){\vector(1,1){60}}
\put(89,158){\vector(-1,-1){60}}

\put(86,11){\vector(-1,1){60}}
\put(30,76){\vector(1,-1){60}}

\end{picture}
}

\put(65,-25)
{

\begin{picture}(100,100)

\put(90,81){${\bf CLis}$}
\put(2,81){${\bf CLis}$}
\put(179,81){${\bf CLis}$}
\put(92,158){${\bf Cde}$}
\put(92,4){${\bf Cde}$}

\put(28,88){\vector(1,0){60}}
\put(88,81){\vector(-1,0){60}}
\put(117,88){\vector(1,0){60}}
\put(177,81){\vector(-1,0){60}}
\put(98,17){\vector(0,1){60}}
\put(104,77){\vector(0,-1){60}}
\put(104,92){\vector(0,1){60}}
\put(98,152){\vector(0,-1){60}}

\put(112,16){\vector(1,1){60}}
\put(176,71){\vector(-1,-1){60}}

\put(177,98){\vector(-1,1){60}}
\put(112,154){\vector(1,-1){60}}

\put(32,93){\vector(1,1){60}}
\put(89,158){\vector(-1,-1){60}}

\put(86,11){\vector(-1,1){60}}
\put(30,76){\vector(1,-1){60}}

\end{picture}
}

\end{picture}

\newpage

\section{Purified Contexts and Doubly Based Lattices}

In order to obtain not only (dual) adjunctions but even categorical (dual)
equivalences between certain categories of contexts and complete lattices,
we must transfer the role played by the object and attribute sets to the
realm of complete lattices. To that aim, we introduce {\em doubly based 
lattices} as triples $({\rm K},J,M)$ where ${\rm K} = (K,\leq )$ is a complete 
lattice, $J$ is a join-dense subset ({\em join-base}) and $M$ is a meet-dense 
subset ({\em meet-base}) of ${\rm K}$. Any context $\KK = (J,M,I)$ gives rise to 
a doubly based lattice
\BM
\B^o\KK = (\B\2 \KK ,J_0,M_0 )  \ $ with $\ J_0 = \ga_{\KK}[J]$ \ and \
$M_0 = \mu_{\KK}[M].
\EM
For any separately continuous morphism 
\BM
(\al ,\be ) : \KK = (J,M,I) \map \LL = (H,N,R)
\EM
the lifted map $\al^{\ra}$ preserves not only joins but also the join-bases,
i.e. $\al^{\ra} [J_0]\inc H_0$, and $\be^{\ra}$ preserves not only meets but
also the meet-bases, i.e. $\be^{\ra} [M_0]\inc N_0$, 
on account of the equations
\BM
\al^{\ra}\!\circ\ga_{\KK} = \ga_{\LL}\circ\al \ $
and $\ \be^{\ra}\!\circ\mu_{\KK} = \mu_{\LL}\circ\be .
\EM
In that way, we obtain a functor $\B^o$ from the category 
category ${\bf Cjm}$ of contexts with separately continuous 
morphisms to the category ${\bf CLjm^o}$ of doubly based lattices with 
mapping pairs $(\vi ,\psi )$ such that $\vi$ preserves joins and the 
selected join-bases, while $\psi$ preserves meets and the selected 
meet-bases.
In the opposite direction, we may assign to each doubly based lattice
${\bf K} = ({\rm K},J,M)$ the {\em base context}
$
\C^{\1 o} {\bf K} = (J,M,\leq )
$
where $\leq$ denotes the given order relation of ${\rm K}$ but also the induced
relation between $J$ and $M$. Then $\C^{\1 o}$ becomes a functor, acting on  
morphisms by restriction to the given join- and meet-bases (see
the proof of Theorem \ref{pur}). A context $\KK$ is said to be 
{\em purified} if $\ga_{\KK}$ and $\mu_{\KK}$ are injective -- in other 
words, if $\eta_{\KK}$ induces an isomorphism between $\KK$ and the context 
$\C^{\1 o} \B^o \KK = (J_0,M_0,\leq )$.
For any doubly based lattice ${\bf K} = ({\rm K},J,M)$, the context 
$\C^{\1 o} {\bf K} = (J,M,\leq )$ is purified. Moreover, by join-density of $J$ and 
meet-density of $M$, we have an isomorphism

\vx

\ \ \ $
\ep_{\bf K} : \B\2  (J,M,\leq )\map {\rm K}\2 , 
\ (A,B)\mapsto\Sup\! A =\Inf B
\EM
with inverse
\BM
\iota_{\bf K} : {\rm K} \map \B\2  (J,M,\leq )\2 ,
\ x\mapsto (J\2\cap \da x , M\2\cap \ua x ) \2 .
\EM
Since $\iota_{\bf K}$ maps $J$ onto $J_0$ and $M$ onto $M_0$, we may regard
that isomorphism as a morphism between ${\bf K}$ and $\B^o \C^{\1 o} {\bf K}$ in 
the category ${\bf CLjm^o}$.

\vx


\BT
\label{pur}
The augmented concept lattice functor $\B^o$ induces an 
equivalence between the category ${\bf Cjm^o}$ of purified contexts with 
separately continuous morphisms and the category ${\bf CLjm^o}$ of doubly 
based lattices with join-meet preserving pairs that induce mappings between the 
respective bases.
In the opposite direction, an equivalence is established by the functor 
$\C^{\1 o}$. 

\vx

\nid Similarly, the category ${\bf C^o}$ of purified contexts 
with conceptual morphisms is equivalent to the category ${\bf CL^o}$
of doubly based lattices with complete homomorphisms preserving the 
join- and meet-bases. 
\ET

\BP
We have already seen that for any ${\bf Cjm^o}$-morphism $(\al ,\be)$, the
map $\B^o (\al ,\be ) = (\al^{\ra},\be^{\ra})$ is a ${\bf CLjm^o}$-morphism.
The equation
\BM
\C^{\1 o}\B^o (\al ,\be )\circ \eta_{\KK} = 
(\al^{\ra}\! \circ \ga_{\KK},\be^{\ra}\!\circ \mu_{\KK})
= (\ga_{\LL}\!\circ\al ,\mu_{\LL}\!\circ\be ) =\eta_{\LL}\circ (\al ,\be ) 
\EM
shows that $\eta$ is a natural
isomorphism between the identity functor on the category of purified contexts and
the composite functor $\C^{\1 o}\B^o$. 

\vx
On the other hand, given an arbitrary ${\bf CLjm^o}$-morphism $\vi$ between doubly
based lattices ${\bf K} =({\rm K},J,M)$ and ${\bf L} = ({\rm L},H,N)$, 
we must show that $\al : J\map H,\ j\mapsto \vi (j)$ and 
$\be : M\map N, \ m\mapsto \vi (m)$ are continuous maps. 
The extents of $\C^{\1 o} {\rm L}$ are of the form $N\cap{\da y}$, and 
\BM
\al^- [N\cap{\da y}] = \{ j\in J : \vi (j)\leq y\} = \{ j\in J : j\leq \vi^*(y)\} = 
J\cap {\da\vi^*(y)}
\EM
is then an extent, too. Thus, $\al$ is extent continuous, and dually, $\be$ is 
intent continuous, so that 
$(\al ,\be )$ is a ${\bf Cjm^o}$-morphism. Moreover, the equation
\BM
(H\2\cap \da \vi (x), N\2\cap \ua \vi (x)) =
(H\2\cap \da\Sup\vi [J\2\cap {\da x}\2 ],N\2\cap\ua\Inf\vi [M\2\cap {\ua x}\2 ])=$

\vx

\ \ \ $
(\vi [J\2\cap\da x]^{\ua\da},\vi [M\2\cap\ua x]^{\da\ua}\2)
\EM
yields the naturality of the isomorphism 
$\iota : {\bf 1}_{\3\bf CLjm^o}\!\map\B^o \C^{\1 o} \2$ and its inverse:

\vx

\ \ \ $
\iota_{\LL}\circ\vi = \B^o \C^{\1 o}\vi \circ \iota_{\KK} \2 , \ $ hence
$\ \vi \circ \ep_{\KK} = \ep_{\LL}\circ\B^o \C^{\1 o} \vi \2 .
\EM
If $\vi$ is a ${\bf CL^o}$-morphism, we have to verify that the mapping pair 
$\C^{\1 o}\vi =(\al , \be)$ from $\C^{\1 o}{\bf K}$ to $\C^{\1 o}{\bf L}$ is conceptual. 
To that aim, we compute for $A\inc J$
and $y = \Sup \vi [A\2 ]$, using join- and meet-density of the bases:
\BM
y = \vi (\Sup A) = \vi (\Inf (M\cap\3 {\ua\!\Sup} A)) = 
\Inf \vi [ M\2\cap {\ua \Sup A}] = \Inf \be [A^{\ua}\2 ]$,

\vx

\ \ \ $\al [A\2 ]^{\ua\da} =(N\cap\ua y)^{\da} =  H\cap {\da y} =
\be [A^{\ua}\2 ]^{\da}.
\EM
A dual reasoning yields the identity 
$\be [B\2 ]^{\da\ua} = \al [B^{\da}\2 ]^{\ua}$ for $B\inc M$.
\EP 

\vx

An obvious question is now whether the above equivalence theorem has an analogue
for concept continuous maps. The answer is in the affirmative, but the choice of
morphisms is a bit subtle. Let us consider the category ${\bf CL_{o*}}$ of doubly 
based lattices with join-base preserving doubly
residuated maps $\psi$ whose double upper adjoint $\psi^{**}$ preserves the 
meet-bases. By passing to these double upper adjoints, one obtains an
isomorphism between ${\bf CL_{o*}}$ and the category ${\bf {CL_o\!}^*}$ of doubly 
based lattices with meet-base preserving doubly residual maps whose double lower 
adjoint preserves the join-bases. A different isomorphism between these
categories results from dualization of the order relations, leaving unchanged
the underlying mappings of the morphisms. But notice that ${\bf CL_{o*}}$ 
and ${\bf {CL_o\!}^*}$ are not dual to the category ${\bf CL^o}$, but to that
category ${\bf CL_o}$ whose morphisms have a join-base preserving lower adjoint 
and a meet-base preserving upper adjoint.

\BT
\label{equicont}
Assigning to each concept continuous morphism $(\al ,\be)$ between purified 
contexts $\KK$ and $\LL$ the map $\al^{\ra} : \B\2\KK \map \B\2\LL$, one obtains 
an equivalence functor $\B_{o*}$ between the category ${\bf C_o}$ of
purified contexts with context continuous pairs and the category ${\bf CL_{o*}}$.
Similarly, sending $(\al ,\be )$ to $\be^{\ra}$, one obtains an equivalence 
functor ${\B_o\!}^*$ between the categories ${\bf C_o}$ and ${\bf {CL_o\!}^*}$.
In the opposite direction, equivalence functors 
from ${\bf CL_{o*}}$ and ${\bf {CL_o\!}^*}$, respectively, to ${\bf C_o}$ are
given by restriction of the doubly residuated maps to the join-bases and of the 
doubly residual maps to the meet-bases. 
\ET

\BP We have seen earlier that for any concept continuous pair $(\al ,\be )$, the 
map $\al^{\ra}$ is doubly residuated, with 
$\al^{\ra *}\! =\al^{\la}\! = \be^{\la}$ and 
$\ \al^{\ra **}\! = \be^{\la *}\! =\be ^{\ra}\2 $,
and that $\al^{\ra}$ preserves the join-bases, while $\be^{\ra}$ 
preserves the meet-bases.

\vx

Given a ${\bf CL_{o*}}$-morphism $\psi = \vi_* : {\bf K}\map {\bf L}$, we show 
that the mapping pair $\C_{o*}\psi =(\al , \be)$ from $\C^{\1 o}{\bf K}$ to 
$\C^{\1 o}{\bf L}$ with $\al = \psi|_J =\vi_*|_J: J\map H$ and 
$\be = \psi^{**\!}|_M =\vi^*\! |_M: M\map N$
is concept continuous. For $C\!\inc\! H$, we have:

\vx

\ \ \ $
g\in \al^-[C^{\ua\da}\2 ] \!\iff\! \al (g)\in H\cap{\da\Sup C} \!\iff\!
\vi_* (g)\leq \Sup C \!\iff\! g\leq \vi (\Sup C)$

\vx

\hspace{-2.5ex}$\iff\alq m\in M\2 (\vi (\Sup C)\leq m \!\imp\! g\leq m)$

\vx

\hspace{-2.5ex}$\iff\alq m\in M\2 (\vi^* (m)\in N\cap{\ua\Sup C}\!\imp\! g\leq m)$

\vx

\hspace{-2.5ex}$\iff\alq m\in M\2 (\be (m)\inc C^{\ua}\!\imp\! g\leq m) \iff 
g\in \be^- [C^{\ua}\2 ]^{\da}.
\EM
This and the dual equation $\be^- [D^{\da\ua}\2 ] = \al^- [D^{\da}\2 ]^{\ua}$
prove concept continuity. 

\vx

\nid That $\eta$ and $\ep$ are natural isomorphisms between the identity functors
and the functors composed by $\C_{o*}$ and $\B_{o*}$ is checked as in the
previous proof. ${\bf {CL_o\!}^*}$-morphisms are treated analogously. 
\EP


\BCO
\label{isoclass}
{\rm (1)} \ Via the augmented concept lattice functor, the isomorphism classes
of purified contexts bijectively correspond to the isomorphism classes of doubly
based lattices. In particular, every purified context is isomorphic to one of the
form $\C^{\1 o}{\bf K}$ for a doubly based lattice ${\bf K}$. 

\vh

{\rm (2)} \ The conceptual morphisms between purified contexts $\KK = \C^{\1 o}{\bf K}$
and $\LL = \C^{\1 o}{\bf L}$ are exactly the restrictions of the base-preserving
complete homomorphisms between the doubly based lattices ${\bf K} \iso \B^o\KK$ 
and ${\bf L}\iso \B^o\LL$.

\vh

{\rm (3)} The concept continuous morphisms between purified contexts 
$\KK = \C^{\1 o}{\bf K}$ and $\LL = \C^{\1 o}{\bf L}$ are exactly the pairs formed by
the base restrictions of the lower adjoint $\vi_*$ and the upper adjoint $\vi^*$
of complete homomorphisms $\vi$ from ${\bf L} \iso \B^o\LL$ 
to ${\bf K}\iso \B^o\KK$ such that $\vi_*$ preserves the join-bases and $\vi^*$
the meet-bases.

\vh

{\rm (4)} \ Conversely, every base-preserving complete homomorphism between
doubly based lattices is induced by a unique conceptual morphism between the 
underlying base contexts, and every complete homomorphism with join-base preserving 
lower adjoint and meet-base preserving upper adjoint is induced by a unique
concept continuous morphism in the opposite direction.
\ECO

\newpage

\section{Dualities and Galois Connections for Contexts}

As expected, there are not only equivalences but also {\em dualities} between 
certain categories of contexts and complete lattices. In most cases, such 
dualities are obtained by composing the already established equivalences with the 
dual isomorphisms resulting from the passage between lower and upper adjoints. 
For example, the category ${\bf CL^o}$ of doubly based lattices with 
base-preserving complete homomorphisms is dually isomorphic to the following two 
categories with the same objects: the morphisms of ${\bf CL^{o*}}$ are those 
doubly residual maps whose lower adjoint preserves join- and meet-bases (hence is 
a ${\bf CL^o}$-morphism), and the morphisms in ${\bf {CL^o\!}_*}$ are those doubly 
residuated maps whose upper adjoint preserves join- and meet-bases.
Thus, Theorem \ref{pur} amounts to:

\BCO
\label{dualconc}
Sending each conceptual morphism $(\al ,\be )$ to the doubly residual map
$\al^{\la}$, one obtains a dual equivalence between the category ${\bf C^{\2 o}}$ 
of purified contexts with conceptual morphisms and the category ${\bf CL^{o*}}$, 
while sending $(\al ,\be )$ to the doubly residuated map $\be^{\la}$, one obtains
a dual equivalence between the categories ${\bf C^o}$ 
and ${\bf {CL^o\!}_*}$.
\ECO

\vx

As already observed earlier, in the same way, the category ${\bf CL_o}$ of doubly
based lattices with complete homomorphisms $\vi$ having a join-base preserving
lower adjoint $\vi_*$ and a meet-base preserving upper adjoint $\vi^*$
is dually isomorphic to the categories ${\bf CL_{o*}}$ and ${\bf {CL_o\!}^*}$. 
Thus, from Theorem \ref{equicont}, we immediately derive:

\BCO
\label{dualcont}
Sending each concept continuous pair $(\al ,\be )$ to the complete
homomorphism $\al^{\la}\! =\be^{\la}$, one obtains a dual equivalence between the 
category ${\bf C_o}$ of purified contexts with concept continuous morphisms and 
the category ${\bf CL_o}$.
\ECO

A basic remark about canonical order structures on contexts is now long overdue.
Both the objects and the attributes of any context carry a natural 
{\em ``specialization order''}, given by
\BM
j\leq k \!\iff\! k^{\ua}\inc j^{\ua}\ \ $ and 
$\ \ m\leq n \!\iff\! m^{\da}\inc n^{\da}.
\EM 
These two relations are obviously quasi-orders (reflexive and transitive), and 
they are partial orders (antisymmetric) iff the context is purified.
The restriction to purified contexts has great structural advantages but causes 
no essential loss of generality, because every context $\KK = (J,M,I)$ has a 
{\em purification} $\C^{\1 o}\B^o\KK = (J_0,M_0,\leq )$ whose concept lattice is 
isomorphic to that of the original context. Note the following implication:
\BM
g\!\leq\! j \3 I \3 m \!\leq\! n \imp g\3 I\3 n .
\EM
It is also convenient to know that for any base context 
$\C^{\1 o}{\bf K} = (J,M,\leq )$, the specialization orders are induced by the lattice
order.

\vx

The following identities connecting
the partners of mapping pairs between purified contexts are easily verified with 
the help of Lemmas \ref{conceptual} and \ref{concont} (maxima and minima refer 
to the specialization orders):

\BPR
Each partner of a conceptual mapping pair $(\al ,\be )$ between purified contexts
determines the other uniquely, by the identities
\BM
\al (g) = \max\3 \{ h : g^{\ua}\inc \be^- [h^{\ua}\2 ]\}\3 , \ \ 
\be (m) = \min\3 \{ n : m^{\da}\inc \al^- [n^{\da}\2 ]\}\3 .
\EM
Similarly, each partner of a concept continuous mapping pair $(\al ,\be )$ between
purified contexts determines the other uniquely, by the identities
\BM
\al (g) = \min\3 \{ h : \be^- [h^{\ua}\2 ]\inc g^{\ua}\}\3 , \ \ 
\be (m) = \max\3 \{ n : \al^- [n^{\da}\2 ]\inc m^{\da}\}\3 .
\EM
\EPR

\nid On account of these facts, it would suffice to consider single maps between 
the object or the attribute sets of purified contexts as conceptual or concept 
continuous morphisms. However, the approach via mapping pairs makes the interplay
between object and attribute sets of contexts more transparent.

\vx

In view of the striking similarities between conceptual and concept continuous
morphisms, and encouraged by the dual isomorphisms between the corresponding 
lattice categories ${\bf CLc}$ and ${\bf CLc\!\2_*^*}$ etc.
one might wish to find similar dualities between suitable subcategories of
the context categories ${\bf C^o}\inc {\bf Cc}$ and 
${\bf C_o}\inc {\bf Cc\!\2_*^*}$. This is in fact possible, 
but the appropriate choice of morphisms might look a bit technical at first glance. 
However, from the Galois-theoretical point of view, it is rather natural.
By slight abuse of language, we call a mapping pair $(\al ,\be )$ between
purified contexts $\KK$ and $\LL$ {\em residuated} if it is concept 
continuous, $\al$ is residuated and $\be$ is residual (i.e. dually residuated)
with respect to the specialization orders. 
On the other hand, we say $(\al ,\be )$ is {\em residual} if it is conceptual, 
$\al$ is residual and $\be$ is residuated (i.e. dually residual). The resulting 
categories of purified contexts are denoted by ${\bf C_r}$ and ${\bf C^r}$, 
respectively. 
As morphisms in the corresponding category ${\bf CLr}$ of doubly based
lattices we take the {\em base-preserving and -reflecting} complete
homomorphisms, i.e. those maps which do not only preserve joins, meets, 
join-bases and meet-bases, but also have the property that their
lower adjoint preserves the join-bases and their upper adjoint
preserves the meet-bases, too. 

\vx

\BL
\label{res}
If $(\al ,\be ) : \KK\map \LL$ is a mapping pair between contexts such that
$\al : J\map H$ has an upper adjoint $\al^*$ and $\be : M\map N$ has a lower 
adjoint $\be_*$, then the following statements are equivalent:
\BIT
\aa $(\al ,\be)$ is concept continuous (residuated).
\bb $(\al^* ,\be_*)$ is conceptual (residual).
\cc $\al (j)^{\ua} = {\be_{*}}^-[j^{\ua}\2 ]$ for all $j\in J$ and
$\be (m)^{\da} = {\al^{*}}^-[m^{\da}\2 ]$ for all $m\in M$. 
\dd $\al^{*}(h)^{\ua} = \be^- [h^{\ua}\2 ]$ for all $\in H$ and
$\3\be_{*} (n)^{\da} = \al^- [n^{\da}\2 ]$ for all $n\in N$.
\EIT
\EL

\BP
(a)$\!\iff\!$(b). On account of Corollary \ref{isoclass}, we may assume that $\KK$
and $\LL$ are the base contexts of two doubly based lattices ${\bf K}$ and 
${\bf L}$, so that $(\al ,\be )$ is concept continuous (hence residuated) iff 
$\al$ is induced by the lower adjoint and $\be$ by the upper adjoint of a common 
base-preserving complete homomorphism $\vi : {\bf L}\map {\bf K}$. 
Since the specialization orders are induced by the lattice orders, both the upper
adjoint $\al^*$ and the lower adjoint $\be_*$ is induced by $\vi$. 
It follows that $\vi$ is a ${\bf CLr}$-morphism. Again by Corollary 
\ref{isoclass}, this is equivalent to saying that $(\al^*,\be_*)$ is conceptual, 
hence residuated.

\vx

For (a)$\!\iff\!$(d), use the identities 
$\al^- [n^{\da}\2 ] = \be^- [n^{\da\ua}\2]^{\da}$ and
$\be^- [h^{\ua}\2 ] = \al^- [h^{\ua\da}\2 ]^{\ua}$ characterizing concept 
continuity.

\vh

\nid That (c) and (d) are equivalent is immediate from the equivalences

\vx

\BTA {cccccc}
$j\in \al^- [n^{\da}\2 ]$&$\!\!\!\iff\!\!\!$&$n\in \al(j)^{\ua}\3 ,$&
$n\in {\be_*}^- [j^{\ua}\2 ]$&$\!\!\!\iff\!\!\!$&$j\in \be_*(n)^{\da}$\\
$h\in \al^{*-} [m^{\da}\2 ]$&$\!\!\!\iff\!\!\!$&$m\in \al^*(h)^{\ua}\3 ,$& 
$m\in \be^- [h^{\ua}\2 ]$&$\!\!\!\iff\!\!\!$&$h\in \be (m)^{\ua}$.\hspace{5ex}\EP
\ETA


\nid From the characterizations (c) and (d) it is obvious that the category of 
posets with residuated and residual maps is embedded in ${\bf C^r}$ (by sending 
$\vi$ to $(\vi ,\vi )$), and the category of posets with doubly residuated (or 
residual) maps is embedded in ${\bf C_r}$ (by sending $\psi$ to $(\psi ,\psi^{**}$)
or $(\psi_{**},\psi )$, respectively). 

\BPR
\label{residuated}
For a mapping pair $(\al ,\be )$ between purified contexts $\KK$ and $\LL$,
the following conditions are equivalent:
\BIT
\aa $(\al ,\be )$ is residuated, that is, a ${\bf C_r}$-morphism.
\IT[${\rm (a^*)}$]  $\al$ has an upper adjoint, $\be$ has a lower adjoint,
and $(\al^*,\be_*)$ is residual, hence a ${\bf C^r}$-morphism.
\bb $(\al ,\be )$ reflects incidence, and there is an incidence
preserving mapping pair $(\al^{\bullet} ,\be_{\bullet}): \LL \map \KK$ with
$\al^{\bullet}(h)^{\ua} = \be^- [h^{\ua}\2 ]$ and
$\be_{\bullet} (n)^{\da} = \al^- [n^{\da}\2 ]$.
\IT[${\rm (b^*)}$] $(\al ,\be )$ reflects incidence, and there is an incidence
preserving mapping pair $(\al^{\bullet} ,\be_{\bullet}): \LL \map \KK$ with
$\al (j)^{\ua} = {\be_{\bullet}}^-[j^{\ua}\2 ]$ and
$\be (m)^{\da} = {\al^{\bullet}}^-[m^{\da}\2 ]$.
\cc There is a unique ${\bf CLr}$-morphism $\vi :\B^o \LL \map \B^o \KK$
such that\\
$\ga_{\LL}\!\circ\al = \vi_*\!\circ\ga_{\KK}$ \ and \ 
$\mu_{\LL}\!\circ\be = \vi^*\!\circ\mu_{\KK}$.
\IT[${\rm (c^*)}$] There is a unique ${\bf CLr}$-morphism $\vi :\B^o \LL \map \B^o \KK$
such that\\
$\vi\circ\ga_{\LL} = \ga_{\KK}\!\circ\al^*$ and 
$\vi\circ\mu_{\LL} = \mu_{\KK}\!\circ\be_*.$
\EIT
Moreover, if these conditions hold then $\al^{\bullet} =\al^*$,
$\be_{\bullet} = \be_*$, and

\vh

\ \ \ $
\vi = \al^{\la\!} = \al^{*\ra\!} = \al^{\ra *\!} = \be^{\la\!} = {\be_*\!}^{\ra}
= {\be^{\ra}\!}_*\3 .
\EM
\EPR

\BPI (300,90)

\put(100,80){$\KK$}
\put(113,86){\vector(1,0){60}}
\put(173,80){\vector(-1,0){60}}

\put(130,90){$(\al ,\be )$}
\put(125,70){$(\al^* ,\be_* )$}
\put(177,80){$\LL$}

\put(104,77){\vector(0,-1){60}}
\put(32,46){$\eta_{\KK} =(\ga_{\KK},\mu_{\KK})$}
\put(75,5){$\C^{\1 o}\B^o\KK$}

\put(179,77){\vector(0,-1){60}}
\put(182,46){$\eta_{\LL} = (\ga_{\LL},\mu_{\LL})$}
\put(177,5){$\C^{\1 o}\B^o\LL$}
\put(173,6){\vector(-1,0){60}}
\put(113,12){\vector(1,0){60}}
\put(132,-4){$(\vi ,\vi )$}
\put(127,16){$(\vi_* ,\vi^* )$}

\EPI

\BP
\nid (a)\im(b). 
From Lemma \ref{res}, we 
infer that $\al^{\bullet}=\al^*$ and $\be_{\bullet}=\be_*$ have
the desired properties; $(\al^*,\be_*)$ preserves 
incidence, because $n\in h^{\ua}$ together with $n\leq\be (\be_* (n))$ implies 
$\be (\be_*(n))\in h^{\ua}$, hence $\be_*(n)\in \be^- [h^{\ua}\2 ]=\al^*(h)^{\ua}$. 
That this is the only possible choice for
$\al^{\bullet}$ and $\be_{\bullet}$ may be checked as follows. First, the
given identities show that $\al$ and $\be$ are continuous, hence isotone
with respect to the specialization orders. But $\al^{\bullet}$ is 
isotone, too: $h\leq k$ means $k^{\ua}\inc h^{\ua}$, which
entails $\al^{\bullet}(k)^{\ua}=\be^-[k^{\ua}\2 ]\inc\be^-[h^{\ua}\2 ]=
\al^{\bullet}(h)^{\ua}$, that is, $\al^{\bullet}(h)\leq \al^{\bullet}(k)$.
Next, we have $\al^{\bullet}(\al (j))^{\ua} = \be^- [\al (j)^{\ua}\2 ]
\inc j^{\ua}$ since $(\al ,\be )$ reflects incidence; therefore, 
$j\leq \al^{\bullet}(\al (j))$. 
On the other hand, the inequality $\al (\al^{\bullet} (h))\leq h$ follows
from the hypothesis that $(\al^{\bullet},\be_{\bullet})$ preserves 
incidence:
$n\in h^{\ua}$ implies $\al^{\bullet}(h)\in \be_{\bullet}(n)^{\da} = 
\al^- [n^{\da}\2 ]$ and then $n\in \al (\al^{\bullet} (h))^{\ua}$. 
Thus, $\al^{\bullet}$ is the upper adjoint of $\al$, and similarly for 
$\be_{\bullet}$.

\vh

\nid (a)$\!\iff\!$(a$^*$) \ and \ (b)$\!\iff\!$(b$^*$) also follow from 
Lemma \ref{res}.

\vh

\nid (b)\im(c) and (c$^*$). 
Lemma \ref{res} tells us that $(\al ,\be )$ is concept continuous.
By Corollary \ref{charcon}, there is a unique complete homomorphism
$\vi =\al^{\la}\! = \be^{\la}$ from $\B\2 \LL$ to $\B\2 \KK$ with 
$\ga_{\LL}\!\circ\al = \vi_*\!\circ\ga_{\KK}$ \ and \ 
$\mu_{\LL}\!\circ\be = \vi^*\!\circ\mu_{\KK}$.
Moreover, the equation
$\vi (\ga_{\LL} (h)) = (...,\be^- [h^{\ua}\2 ]) = (...,\al^* (h)^{\ua}) =
\ga_{\KK} (\al^*(h))$
and its dual yield the identities
$\vi\circ\ga_{\LL} = \ga_{\KK}\!\circ\al^*$  and  
$\vi\circ\mu_{\LL} = \mu_{\KK}\!\circ\be_*.$
In particular, $\vi$ preserves the join- and meet-bases. 

\vh

\nid (c)\im(a). By Corollary \ref{charcon}, the pair $(\al ,\be )$ is
concept continuous. For $h\in H$, we have
$(\al^- [h^{\ua\da}\2 ], \al^- [h^{\ua\da}\2 ]^{\ua}) =
\al^{\la} (h^{\ua\da},h^{\ua}) = \vi(\ga_{\LL}(h))\in \ga_{\KK}[J],$
hence 
$\al^- [\da h] = \al^- [h^{\ua\da}\2 ] = j^{\ua\da} = {\da j}$ for some 
$j\in J$; thus, $\al$ is residuated. Dually, one shows that $\be$ is 
residual.

\vh

\nid (c$^*$)\im(a$^*$). Use Corollary \ref{char} (with $\KK$ and $\LL$ exchanged).

\vh



\vh

\nid (b$^*$)\im(c$^*$) is established in the same manner as (b)\im(c).
\EP

\vx

\nid Let us note a few additional properties of residuated mapping pairs.

\BCO
If $(\al ,\be)$ is a residuated pair then $(\al^* ,\be_*)$ is a residual
pair, and the following identities are fulfilled: 
\BM
\al^{*-} [A^{\ua\da}\2 ] = \be [A^{\ua}\3 ]^{\da}\3 , \ \ \
{\be_*}^-[B^{\da\ua}\2 ] = \al [B^{\da}\2 ]^{\ua},$

\vx

\ \ \ $\al^* [C\2 ]^{\ua\da}\2 = \al^- [C^{\1\ua\da}\3 ] = 
\be_* [C^{\ua}\2 ]^{\da} = \be^- [C^{\ua}\2 ]^{\da},$

\vx

\ \ \ $\al^* [D^{\da}\2 ]^ {\ua} = \al^- [D^{\da}\2 ]^{\ua} =
\be_* [D]^{\da\ua} = \be^- [D^{\da\ua}\2 ]\2 .
\EM
\ECO

\BP
For the first equation, observe the equivalences
\BM
h\in \al^{*-} [A^{\ua\da}\2 ] \!\iff\! A^{\ua}\inc\al^*(h)^{\ua}=
\be^- [h^{\ua}\2 ] \!\iff\! \be [A^{\ua}\2 ]\inc h^{\ua} \!\iff\!
h\in\be [A^{\ua}\2 ]^{\da}.
\EM
The identity $\al^* [C]^{\ua\da} = \be_* [C^{\ua}\2 ]^{\da}$ follows 
from conceptuality of $(\al^* ,\be_* )$ (see Lemma \ref{conceptual}), and 
the identity $\al^- [C^{\ua\da}\2 ] = \be^- [C^{\ua}\2 ]^{\da}$ from concept
continuity of $(\al ,\be )$ (see Lemma \ref{concont}). For
$\al^- [C^{\ua\da}\2 ] = \be_* [C^{\ua}\2 ]^{\da}$, use the equivalences
\BM
j\in \al^- [C^{\ua\da}\2 ] \!\iff\! C^{\ua}\inc \al (j)^{\ua}
= {\be_*}^-[j^{\ua}\2 ]\!\iff\! \be_* [C^{\ua}\2 ] \inc j^{\ua}
\!\iff\! j\in \be_* [C^{\ua}\2 ]^{\da}\2 .
\EM
The othe equations are derived analogously.
\EP

\vx

\nid From the equivalence of (a) and (a$^*$) in Proposition \ref{residuated},
we conclude:

\BCO
\label{dualr}
By passing to adjoints, the context categories ${\bf C_r}$ and  ${\bf C^r}$ 
are dually isomorphic to each other.
\ECO

Now, we are in a position to establish the main result of this section:

\BT
\label{dual}
Assigning to each residual pair $(\al ,\be )$ the complete
homomorphism $\al^{\ra}\! = \be^{\ra}$, one obtains
an equivalence $\B^{\2 r}$ between the category ${\bf C^r}$ of purified 
contexts and the category ${\bf CLr}$ of doubly based lattices. In the opposite
direction, the equivalence functor $\C^{\2 r}$ sends any ${\bf CLr}$-morphism
$\vi$ to the mapping pair built by the restrictions 
of $\vi$ to the join- and meet-bases.

\vh

Similarly, associating with any residuated pair $(\al ,\be )$ the complete
homomorphism $\al^{\la}\! = \be^{\la}$, one obtains
a dual equivalence $\B_{\1 r}^{\2 *}$ between the categories ${\bf C_r}$ and 
${\bf CLr}$. In the opposite direction, the dual equivalence functor 
$\C_{\1 r}^{\2 *}$ sends any ${\bf CLr}$-morphism $\vi$ to the pair 
constituted by the restriction of $\vi_*$ to the join-bases and the restriction of 
$\vi^*$ to the meet-bases.
\ET

\BP
We already know that for any context $\KK = (J,M,I)\2$, the triple
$\B^{\2 r}\KK = \B^{\2 *}_{\1 r}\KK = \B^o\KK =(\B\2 \KK, J_0,M_0)$ is a doubly 
based lattice, and that
\BM
\eta_{\KK} = (\ga_{\KK},\mu_{\KK}) :\KK \map \C^{\1 o}\B^o\KK , 
\ x\mapsto (J\cap {\da x},M\cap {\ua x})
\EM
is a natural isomorphism provided $\KK$ is purified.
On the other hand, for an arbitrary doubly based lattice ${\bf K}=({\rm K},J,M)$,
we have the purified context $\C^{\2 r}{\bf K} = 
 \C_{\1 r}^{\2 *}{\bf K} = \C^{\1 o}{\bf K} =(J,M,\leq )$ 
and the natural isomorphism 
\BM
\iota_{\bf K} : {\bf K}\map \B^o\C^{\1 o}{\bf K}\2 , 
\ x\mapsto (J\cap {\da x},M\cap {\ua x})\2 .
\EM
That for any ${\bf C_r}$-morphism $(\al ,\be )$ the pair 
$\B^{\2 r} (\al ,\be ) = \al^{\la} = \be^{\la}$ has the required 
properties of a ${\bf CLr}$-morphism was shown in Proposition \ref{residuated},
and similarly, for any ${\bf C^r}$-morphism $(\al ,\be )$, the pair 
$\B_{\2 r} (\al ,\be ) = \al^{\ra} = \be^{\ra}$ is a ${\bf CLr}$-morphism, too.

\vx

Conversely, given any ${\bf CLr}$-morphism $\vi : {\bf L}\map {\bf K}$
between doubly based lattices ${\bf L} =({\rm L},H,N)$ and ${\bf K} = 
({\rm K},J,M)$, we have that the restricted maps $\al = \vi_*: J\map H$ and
$\be =\vi^*: M\map N$ form a residuated pair $\C_r^{\2 *}\vi =(\al ,\be)$, 
hence a ${\bf C_r}$-morphism. 
This and the remaining statements are easy consequences of earlier results.
\EP

Let us finally put together all pieces of the Galois duality puzzle.
In the diagram on the next page, we place 13 different categories in three
triangular levels; all categories of one level are mutually equivalent or dual.
Each double line symbolizes a categorical equivalence, while each (non-dotted)
single line stands for a duality. In the table of morphisms, 

\vh

$\vdash \3 J\4$ indicates that join-bases are preserved,

$\vdash M$ indicates that meet-bases are preserved,

$\vi_*$ denotes the lower adjoint and $\vi^*$ the upper adjoint of $\vi$.

\newpage
\BC
{\bf Equivalent and dual categories of contexts and complete lattices}
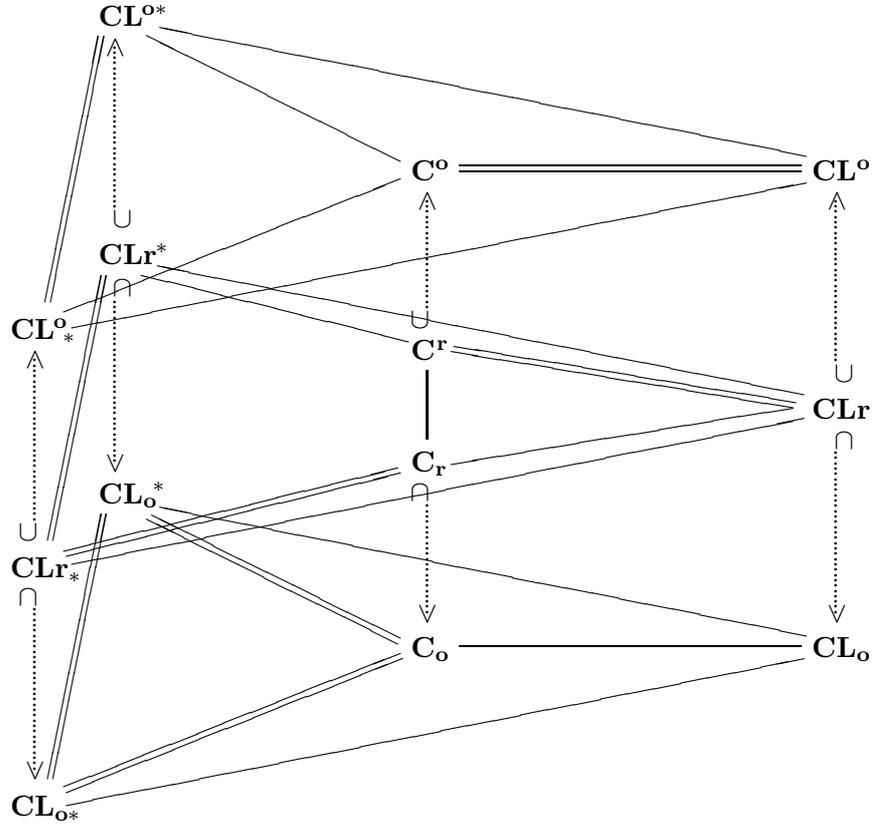
\EC

\BPI (300,320)(0,0)

\put(150,240){${\bf C^o}$}
\put(0,180){${\bf CL^{\!\2 o}\!_*}$}
\put(33,298){${\bf {CL\!\3}^{o*}}$}
\put(150,172){${\bf C^{\1 r}}$}
\put(300,240){${\bf CL^o}$}

\put(0,90){${\bf CLr_*}$}
\put(33,208){${\bf {CLr\!\3}^*}$}
\put(150,130){${\bf C_r}$}
\put(300,150){${\bf CLr}$}

\put(0,0){${\bf CL_{o*}}$}
\put(33,118){${\bf {CL\!\4_o\!}^*}$}
\put(150,60){${\bf C_o}$}
\put(300,60){${\bf CL_o}$}

\put(3,79.5){$\cap$}
\multiput(9,18)(0,2){31}{\circle*{.5}}
\put(5.7,14){$\ve$}

\put(3,104){$\cup$}
\multiput(9,112)(0,2){31}{\circle*{.5}}
\put(5.7,169){$\we$}

\put(308.4,164){$\cup$}
\multiput(309,172)(0,2){31}{\circle*{.5}}
\put(305.7,229){$\we$}

\put(308.4,138.5){$\cap$}
\multiput(309,77)(0,2){31}{\circle*{.5}}
\put(305.7,73){$\ve$}

\put(38.3,222){$\cup$}
\multiput(39,230)(0,2){31}{\circle*{.5}}
\put(35.7,287){$\we$}

\put(38.3,196){$\cap$}
\multiput(39,134)(0,2){31}{\circle*{.5}}
\put(35.7,130){$\ve$}

\put(149.7,184){$\cup$}
\multiput(156,192)(0,2){21}{\circle*{.5}}
\put(152.5,229){$\we$}

\put(149.7,119){$\cap$}
\multiput(156,77)(0,2){21}{\circle*{.5}}
\put(152.5,73){$\ve$}

\put(21,4){\line(5,1){276}}
\put(168,64){\line(1,0){128}}
\put(298,68){\line(-5,1){242}}
\put(14,14){\line(1,5){20}}
\put(16,14){\line(1,5){20}}

\put(168,243){\line(1,0){128}}
\put(168,245.5){\line(1,0){128}}

\put(149,179){\line(-4,1){100}}

\put(294,156){\line(-6,1){128}}
\put(293,154){\line(-6,1){128}}
\put(293,154){\line(-6,-1){128}}

\put(23,95){\line(5,1){274}}
\put(156,142){\line(0,1){26}}
\put(297,159){\line(-5,1){242}}
\put(14,104){\line(1,5){20}}
\put(16,104){\line(1,5){20}}

\put(23,184){\line(5,1){275}}
\put(298,249){\line(-5,1){242}}
\put(13,194){\line(1,5){20}}
\put(15,194){\line(1,5){20}}

\put(20,190){\line(5,2){126}}

\put(20,11){\line(5,2){126}}
\put(21,9){\line(5,2){126}}

\put(20,100){\line(4,1){125}}
\put(21,98){\line(4,1){125}}

\put(146,247){\line(-2,1){95}}

\put(146,67){\line(-2,1){93}}
\put(145,65){\line(-2,1){94}}

\EPI

\vspace{3ex}

\BTA {|c|c|c|c|c|}
\hline
$\hspace{-.5ex}${\em category}$\hspace{-.5ex}$ & {\em objects} & {\em morphisms}
& \multicolumn{2}{c|}{\em additional properties} \\
\hline
\hline
${\bf C^o}$      & purified       & conceptual  &          & \\
\cline{1-1}
\cline{5-5}
${\bf C^{\2 r}}$ & or             & pairs       &          & residual\\
\cline{1-1}
\cline{3-5}
${\bf C_o}$      & reduced        & concept     &          & \\
\cline{1-1}
\cline{5-5} 
${\bf C_{\1 r}}$ & contexts  & continuous pairs &          & residuated\\
\hline
\hline
                 &                &             &$\vdash J\ \ \ \vdash M$     & 
$\vdash J \ \ \ \vdash M$\\
\hline
${\bf CL^ {o*}}$ &                & $\psi$      & $\psi_* \ \ \ \ \ \ \psi_*$ &\\
\cline{1-1}
\cline{4-5}
$\2{\bf CLr^*}$    & doubly       & doubly      & $\psi_* \ \ \ \ \ \ \psi_*$ & 
$\!\!\!\!\psi_{**}\ \ \ \ \ \psi$\\
\cline{1-1}
\cline{4-5}
$\5{\bf {CL_o}^*}\3$ & based      & residual    &                             & 
$\!\!\!\!\psi_{**}\ \ \ \ \ \psi$\\
\cline{1-1}
\cline{3-5}
$\hspace{-1ex}\2{\bf CL^{o\1}}$&or&$\vi$        & $\vi \ \ \ \ \ \ \ \vi$     &
\\
\cline{1-1}
\cline{4-5}
${\!\! \bf CLr}$ & irreducibly    & complete    & $\vi$ \ \ \ \ \ \ $\vi$     & 
$\!\vi_* \ \ \ \ \ \ \vi^*$\\
\cline{1-1}
\cline{4-5}
$\!\!{\bf CL_{o}}$&bigenerated    &homomorphisms&                             &
$\!\vi_* \ \ \ \ \ \ \vi^*$\\
\cline{1-1}
\cline{3-5}
${\bf {CL^o\!\2}_*}$&lattices     & $\psi$      & $\psi^* \ \ \ \ \ \psi^*$   &\\
\cline{1-1}
\cline{4-5}
${\bf CLr\!\2_*}$       &         & doubly      & $\psi^* \ \ \ \ \ \psi^*$   &
$\psi \ \ \ \ \ \ \psi^{**}$\\
\cline{1-1}
\cline{4-5}
${\bf CL_{o *}}$        &         & residuated  &                             &
$\psi \ \ \ \ \ \ \psi^{**}$\\
\hline
\ETA

\vspace{2ex}

\newpage

\section{Reduced Contexts and Irreducibly Bigenerated Lattices}

\vx

We have seen that join- and meet-bases play a crucial role in the passage between
context and complete lattices.
The situation is simplified considerably if we focus on {\em irreducibly 
bigenerated} lattices; these have a least join-base $\J (L)$, consisting 
of all join-irreducibles, and a least meet-base $\M (L)$, consisting of all 
meet-irreducibles. Of course, all finite lattices have that property. 
Any irreducibly bigenerated lattice is isomorphic to the concept lattice of
an up to isomorphism unique reduced context, viz. the {\em standard
context} $\S L = (\J (L),\M (L),\leq)$. An arbitrary context $\KK =(J,M,I)$
is {\em reduced} iff it is purified, each object concept $\ga (j)$ is 
join-irreducible, and each attribute concept $\mu (m)$ is meet-irreducible
in the concept lattice $\B\2 \KK$ -- in other words, iff $\eta_{\KK}$ 
induces an isomorphism between the contexts $\KK$ and $\S\B\2 \KK$.
On the other hand, a complete lattice is irreducibly bigenerated iff it is
isomorphic to the concept lattice of its standard context.
We may regard $\S$ as a covariant functor, sending each join- and meet-irreducibility
preserving complete homomorphism to the pair of its restrictions to the
least join- and meet-bases, respectively. But, of course, there is also a
contravariant standard context functor, restricting any complete 
homomorphism whose lower adjoint preserves join-irreducibility and whose upper 
adjoint preserves meet-ireducibility, to the respective least bases.
Now, the equivalences and dualities between categories of purified contexts and 
doubly based lattices derived in the previous sections immediately lead to the 
following more restricted but technically simpler results:

\BT
\label{redcon}
Under the concept lattice functor and the standard context functor
in the reverse direction, the category of reduced contexts and 
conceptual morphisms is equivalent to the category of irreducibly 
bigenerated lattices and complete homomorphisms preserving the least 
join- and meet-bases. 
Hence, the conceptual morphisms between reduced contexts are in one-to-one
correspondence with the irreducibility preserving complete homomorphisms 
between their concept lattices. 
\ET

\BT
Via the contravariant concept lattice functor and the contravariant
standard context functor, the category of reduced contexts with 
concept continuous pairs is dual to the category of irreducibly bigenerated
lattices and complete homomorphisms whose lower adjoint preserves the join-bases 
and whose upper adjoint preserves the meet-bases. 

Similarly, the category of reduced contexts with residual (respectively 
residuated) mapping pairs is equivalent (respectively dual) to the category of 
irreducibly bigenerated lattices with complete homomorphisms preserving and 
reflecting the least bases.
\ET

\end{document}